\documentclass{article}

\usepackage{arxiv}

\usepackage[utf8]{inputenc} 
\usepackage[T1]{fontenc}    
\usepackage{hyperref}       
\usepackage{url}            
\usepackage{booktabs}       
\usepackage{amsfonts}       
\usepackage{nicefrac}   
\usepackage{microtype}      
\usepackage{lipsum}
\usepackage{graphicx}
\usepackage{subfiles,physics,amsmath,amsthm,amssymb,mathtools,xcolor,bm,url,tikz,tikz-cd,multicol,graphicx,wrapfig,subcaption}

\usepackage{algorithm}
\usepackage{algpseudocode}

\newcommand{\R}{\mathbb{R}}
\newcommand{\C}{\mathbb{C}}

\newcommand{\OO}{\mathcal{O}}
\newcommand{\Sbb}{\mathbb{S}}
\newcommand{\xx}{\pmb{x}}
\newcommand{\rpmb}{\pmb{r}}
\newcommand{\spmb}{\pmb{s}}
\newcommand{\nor}{\text{nor}}

\newcommand{\defeq}{\mathrel{\mathop:}=}

\newcommand{\inc}{\mathrm{in}}
\newcommand{\sca}{\mathrm{sc}}

\newcommand{\sfx}{\mathsf{x}}
\newcommand{\sfF}{\mathsf{F}}
\newcommand{\sfK}{\mathsf{K}}
\newcommand{\sfC}{\mathsf{C}}
\newcommand{\sfS}{\mathsf{S}}
\newcommand{\sfG}{\mathsf{G}}
\newcommand{\sfM}{\mathsf{M}}
\newcommand{\sfO}{\mathsf{O}}
\newcommand{\sfH}{\mathsf{H}}
\newcommand{\sfU}{\mathsf{U}}
\newcommand{\sfV}{\mathsf{V}}

\newcommand{\del}{\partial}
\newcommand{\eps}{\epsilon}
\newcommand{\barr}{\overline}

\DeclarePairedDelimiter\inner{\langle}{\rangle}

\graphicspath{ {./images/} }

\title{Solving the wide-band inverse scattering problem via equivariant neural networks}

\author{
 Borong Zhang \\
  Department of Mathematics\\
  University of Wisconsin-Madison\\
  Madison, WI 53706 \\
  \texttt{bzhang388@wisc.edu} \\
   \And
  Leonardo Zepeda-Nunez \\    
  Department of Mathematics\\
  University of Wisconsin-Madison\\
  Madison, WI 53706 \\
Google Research \\
Mountain View, CA 94043\\
\texttt{lzepedanunez@google.com} 
   \And
  Qin Li \\
  Department of Mathematics\\
  University of Wisconsin-Madison\\
  Madison, WI 53706 \\
  \texttt{qinli@math.wisc.edu} 
}

\begin{document}
\maketitle
\begin{abstract}
This paper introduces a novel deep neural network architecture for solving the inverse scattering problem in frequency domain with wide-band data, by directly approximating the inverse map, thus avoiding the expensive optimization loop of classical methods. The architecture is motivated by the filtered back-projection formula in the full aperture regime and with homogeneous background, and it leverages the underlying equivariance of the problem and compressibility of the integral operator. This drastically reduces the number of training parameters, and therefore the computational and sample complexity of the method. In particular, we obtain an architecture whose number of parameters scales sub-linearly with respect to the dimension of the inputs, while its inference complexity scales super-linearly but with very small constants. We provide several numerical tests that show that the current approach results in better reconstruction than optimization-based techniques such as full-waveform inversion, but at a fraction of the cost while being competitive with state-of-the-art machine learning methods.  
\end{abstract}


\section{Introduction}

Inverse wave scattering is a classical inverse problem that utilizes the scattered wave from an impinging probing signal (or \textit{data}) to infer the acoustic properties of the object being probed (or the \textit{to-be-be-reconstructed media}). The problem finds wide applications in radar imaging \cite{radar}, sonar imaging \cite{sonar}, seismic exploration \cite{seismic}, geophysics exploration \cite{geophysics}, {bio-medical} imaging \cite{medical} and so on.

At the theoretical level, the well studied single-frequency (or monochromatic) inverse wave scattering problem has been proved to produce a unique reconstruction of the media \cite{Kirsch}: if the full incoming-to-scattered wave map is given for a fixed frequency, the data can uniquely determine the underlying object, when the media is assumed to be smooth. Unfortunately, this mathematical statement finds little practical {use in the design of inversion algorithms}. The inversion problem, when studied in the theoretical setting requires the knowledge of the full map, while the numerical implementation only exploits a data set of 
finite, {albeit very large in size}. In addition, it is known that inverse scattering problem is severely ill-posed and the quality of the reconstruction \textit{crucially} depends on the quality of the collected data, particularly, when the data is finite. 
Thus, in summary, the inverse scattering problem in the numerical setting is often confronted with two issues: its sample complexity, i.e., the amount of the data required, and the lack of stability of the reconstruction.

%

{The sample complexity is usually related to the parametrization of the medium to be reconstructed. While the mathematical PDE statement guarantees the reconstruction of a function, numerically one needs to represent this function using a finite, but often very large, number of parameters, particularly when the media is expected to have fine-grained structures.

This calls for a large amount of data for the reconstruction due to so-called diffraction limit \cite{garnier2016passive}. In fact, to capture such small features, one needs to probe the medium with higher-frequency waves, with an increasingly larger number of directions, while the scattered waves need to be sampled at an increasingly finer rate. Qualitatively, in two dimensions the amount of data needed scales proportionally to the square of the inverse of the characteristic length, the length that represents the smallest feature in the media one seeks to recover.}

The lack of stability has been mathematically proved~\cite{Hahner} for monochromatic data: a small perturbation in the collected data leads to disastrous inaccuracies in the reconstruction of the media~\cite{ChenDingLiZepeda}. To overcome such a difficulty, within the classical computational pipelines of optimization-based reconstruction techniques, one typically uses wide-band data, i.e., using different frequencies for the probing waves. Unfortunately this brings two other issues, first, the cost of simulating the wave propagation is prohibitive, particularly at high-frequency, due to combination of the Shannon-Nyquist sampling criterion~\cite{whittaker_1915,shanon1948:criterion} and a strict Courant–Friedrichs–Lewy (CFL) condition~\cite{Courant1928berDP}, and second, the problem is still highly non-convex, so non-physical local minima are abundant \cite{Symes_Chen_Minkoff:2020}. During the past decades, multiple numerical strategies have been taken to ease these difficulties, including fast PDE solvers \cite{ZepedaDemanet:the_method_of_polarized_traces,EngquistYing:Sweeping_PML} that seek to alleviate the computational costs, and new optimization pipelines such as recursive linearization \cite{linearrecursion} and full wave inversion \cite{fwi} that seek to avoid some of the spurious local minima in the optimization loop by a hierarchical processing of the data.

{In a nutshell, in order to reconstruct fine-grained details one needs higher frequency data; however, as the frequency increases, the size of the data becomes larger thus rendering each step in the optimization procedure increasingly expensive, and the optimization loss becomes increasingly non-convex.}

Instead of using an iterative method to reconstruct the media, one alternative is to \textit{directly} approximate the map between the to-be-reconstructed media and the data, which would allows us to circumvent some of the issues mentioned above. 
Thus driven by the empirical success of machine learning in approximating high-dimensional, highly non-linear maps in a myriad of applications~\cite{vision, language, speech}, we investigate whether it is possible to \textit{efficiently} approximate this map by leveraging the power of deep learning, when the to-be-reconstructed medium is parametrized by an \textit{unknown} but relatively low-dimensional manifold.

In our setting, the {question we seek to answer in this paper} is:
\begin{center}
\emph{Are neural networks capable of efficiently approximate the map between media and data?}
\end{center}
The answer is high-likely to be positive as neural operator networks~\cite{fno,deeponet} have been developed by a growing number of researchers {reaching different level of success}. Unfortunately, both difficulties, i.e., lack of stability and sample complexity are also present when approximating the map directly, albeit they are expressed differently. The lack of stability induces a severe lack of convexity in the energy landscape of typical optimization losses, rendering the training of the network challenging, while the number of parameters for the map is often very large and increases rapidly as the frequency increases, which in return requires an increasingly larger amount of training data.

Under the theme of neural approximation of the high-dimensional maps, we seek to alleviate the complexity issues mentioned above for the inverse scattering problem. We study if one can build certain properties of media-data relation into the design of a neural architecture, so to enforce these properties and therefore reduce the computational cost, resulting in an efficient method, particularly when using wide-band data, which has been shown to stabilize the training~\cite{MLZ}. Specifically, we look into the equivariant property and butterfly structure of the media-data relation, and we design a neural network that has these features built-in. Such architecture not only respects the underlying physics, interpreted through the PDE (Helmholtz equation), but also helps to reduce the overall computational cost. 

The main contribution of this work is to construct a model that is able to reconstruct relatively fine-grained features accurately such that:
\begin{itemize}
    \item it has {asymptotically} fewer training parameters, which is achieved by leveraging a detailed description of the linearized operator, the underlying equivariance of the problem, and a compression technique tailored for oscillatory phenomena; and
    \item empirically, it requires much smaller training datasets to achieve high-validation accuracy,
\end{itemize} 
particularly when compared to existing approaches \cite{Khoo_YingSwitchNet:2019,FanYing:scattering,MLZ}. 

There have been several new developments for leveraging ML techniques for inverse problems. In \cite{PINN_Inverse_Problems} the authors used the recently introduced paradigm of physics informed neural networks (PINN) to solve for inverse problems in optics. Aggarwal et al.~introduce a model-based image reconstruction framework \cite{MoDL} for MRI reconstruction. The formulation contains a novel data-consistency step that performs conjugate gradient iterations inside the unrolled algorithm. Gilton et al.~proposed in \cite{Neumann_Networks} a novel network based on Neumann expansion series coupled with a hand-crafted pre-conditioner for linear inverse problems, which recast an unrolled algorithm as elements in the Neumann series. In \cite{Mao:2016} Mao et al.~use a deep encoder-decoder network reminiscent of U-nets \cite{U-Net} for image de-noising, using symmetric skip connections. Networks based on the scattering transform has been proposed \cite{Invariant_Scattering_CNN} to take in account translational equivariance in images. In \cite{Framelets} the authors introduced another framework based on frames for inverse problems, which was applied to computer tomography de-noising \cite{Framelets_denosing}. These works pioneered the application of deep learning and explored the possibility of integrating physics of inverse scattering with novel machine learning algorithms.

In this paper we specifically focus in two features of the inverse scattering problem: the complementary low-rank property of the oscillatory integrals in the inversion formula, and the exploitation of equivariance. We point out that these features have been studied in literature although only \textit{separately}. On the algebraical level, constructions of structures such as H-matrices \cite{MNNH2}, Fast Fourier Transforms \cite{FourierNeuralOp,LRC2020} or butterfly factorizations \cite{linear_butterfly,Khoo_YingSwitchNet:2019,Butterfly-Net2,Yingzhou2018,MLZ} have been extensively studied for reducing the computational cost. We are interested in mimicking these constructions in the design of the network. In addition, in \cite{FanYing:scattering}, the authors proposed a network that is rotationally equivariant, which allows them to reduce the number of degrees of freedom, albeit only using monochromatic data.

Our current approach exploits both traits: our design of the  neural architecture relies on equivariant formulation of the linearized operator, which we further compress it using a Butterfly-type structure \cite{Butterfly-factorization:Liu_Xing2020}, while processing several frequencies simultaneously.  

We organize our paper as follows. In Section~\ref{sec:preliminary}, we present the problem setup and explain the details of the media-data relation. We introduce the equivariant property and the butterfly structure of the media-data relation inherited from the PDE. The embedding of these properties into the design of NN architectures will be discussed in Section~\ref{sec:equivariance} and~\ref{sec:butterfly} respectively. The final section contains several numerical experiments showcasing the properties of the proposed method.

\section{Problem setup}\label{sec:preliminary}
The inverse scattering problem, at large, is the inverse problem associated with the Helmholtz equation, one of the simplest models for time-harmonic wave propagation. Despite its simplicity, it retains the core difficulties of models with more complex physics such as seismic or electromagnetic waves, thus making it an ideal test bed for new algorithms. The Helmholtz equation is the Fourier-in-time transform of the constant-density acoustic wave equation and it reads:
\begin{equation}\label{eqn:helmholtz}
\Delta u(\bm{x})+\omega^2n(\xx)u(\bm{x})=0\,,\quad \xx\in\Omega\subset\mathbb{R}^2
\end{equation}
where $u$ is the total wave field, $\omega$ is the probing wave frequency, and $n$ is the refractive index that encodes the media's wave speed. The domain of interest is denoted by $\Omega \subset \mathbb{R}^2$. We assume that the background media is homogeneous and equal to one, i.e., $n(\xx)=1$ for $\xx\notin\Omega$. For simplicity we denote the perturbation  $\eta(\xx)=n(\xx)-1$, thus we have $\text{supp}(\eta(\xx))\subset\Omega$. 

For a given $n$ (or equivalently $\eta$), the forward problem amounts to solving for the scattered wave field under the Sommerfeld radiation condition when the media is impinged by a probing wave. Due to the homogeneous background assumption we consider probing waves as plane waves:
\begin{equation}
u^\text{in}=e^{i\omega \bm{s}\cdot \bm{x}},
\end{equation}
where $\bm{s}\in\mathbb{S}^1$ is a unitary vector denoting the direction of the incoming plane wave. When injected into \eqref{eqn:helmholtz}, it triggers the scattered wave field $u^\text{sc}$. Define
\begin{equation*}
u^\text{sc}= u-u^\text{in}\,,
\end{equation*}
and utilizing the fact that $u^\text{in}$ solves the Helmholtz equation with trivial media $n(\xx)=1$, it is straightforward to show that $u^\sca(\xx;\spmb)$ solves:
\begin{equation}\label{eqn:scattereqn}
\begin{cases}
\Delta u^{\sca}(\bm{x})+\omega^2(1+\eta(\xx)) u^{\sca}(\bm{x})=-\omega^2\eta(\xx)u^\inc, \\
\frac{\del u^{\sca}}{\del \abs{\bm{x}}}-i\omega u^{\sca}=\OO(\abs{\bm{x}}^{-3/2})
    \text{ uniformly in } \frac{\bm{x}}{\abs{\bm{x}}}\in S^1\text{ as }\abs{\bm{x}}\to\infty.
\end{cases}
\end{equation}
The second row in the equation is termed the Sommerfeld radiation condition. It is imposed at the infinity to ensure the well-possedness of the equation~\eqref{eqn:scattereqn}, which further guarantees the uniqueness of the solution to~\eqref{eqn:helmholtz}.
\begin{figure}[H]
\center
\includegraphics[width=0.4\textwidth]{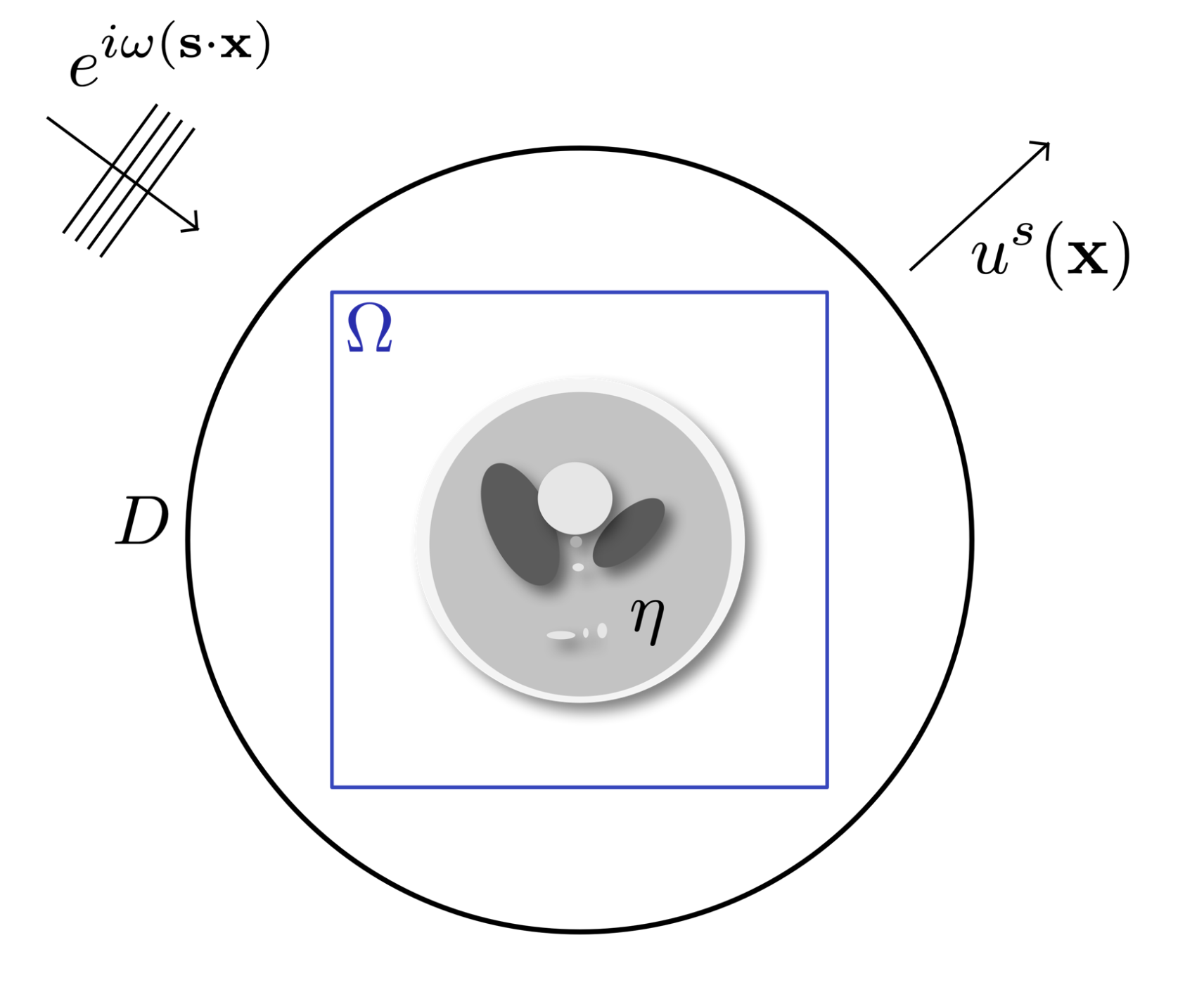}
\caption{The setup for the inverse scattering problem. In the illustration, the media $\eta$, in the domain of interest $\Omega$, is impinged by the probing wave with frequency $\omega$ from the direction \textbf{s}. The scattered field $u^{\textbf{s}}(\mathbf{x})$ is sampled on the disk $D$.}
\end{figure}

The measurements are taken at the circle $R$ away from the origin, i.e., we measure the value of $u^\sca(R\bm{r})$ with $R>\text{radius}(\Omega)$ and $\bm{r}\in\mathbb{S}^1$. We define the continuous far-field pattern (the data) as $\Lambda^\omega(\spmb,\rpmb)$, namely:
\begin{equation} \label{eq:far_field_pattern}
    \Lambda^\omega(\spmb,\rpmb)= u^\sca(R\rpmb;\spmb)\,.
\end{equation}

Since this data is uniquely determined by the configuration of $\eta(\xx)$, we define the map from media $\eta$ to the data $\Lambda^\omega$ as 
\[
\Lambda^\omega = \mathcal{F}^\omega[\eta]\,.
\]
{
\subsubsection*{Problem Formulation}
While the forward problem is to compute $\Lambda^\omega$ for the given $\eta$ (or equivalently to apply $\mathcal{F}$), the inverse problem is to infer $\eta$ from the measured data $\Lambda^\omega$ (or equivalently, to invert $\mathcal{F}$ for $\eta$), which we write formally as
\begin{equation} \label{eq:inv_formulation}
\eta^\ast = \mathcal{F}^{-1}(\{\Lambda^{\omega}\}_{\omega \in \bar{\Omega}}),
\end{equation}
where $\bar{\Omega}$ is a discrete set of frequencies chosen before hand. We note that even though the PDE itself is linear, $\mathcal{F}^\omega$ is a nonlinear map.}

One classical way to numerically execute the inversion is to recast this problem as a PDE-constrained optimization~\cite{optimization}. It seeks the configuration of $\eta$ so that the synthetic data generated by $\eta$ (by solving the PDE~\eqref{eqn:scattereqn}) minimizes the mismatch with respect to the data $\Lambda^\omega$ generated by the probing wave of frequency $\omega$. Such procedure is formulated as
\begin{equation} \label{eq:inv_problem}
\eta^\ast = \mathcal{F}^{-1}(\{\Lambda^{\omega}\}_{\omega \in \bar{\Omega}}) := \text{argmin}_{\nu} \sum_{\omega \in \bar{\Omega}} \| \mathcal{F}^{\omega}[\nu] - \Lambda^{\omega} \|^2.
\end{equation}
The problem is usually solved using highly tailored gradient-based optimization techniques, where the gradient is computed using adjoint-state methods \cite{Plessix_2006:ajoint_state}.

\subsection{Filtered Back-Projection}

Instead of using optimization to solve \eqref{eq:inv_problem} we seek to approximate the map directly using a neural network. As a motivation for our architecture, we briefly introduce a classical strategy that hinges on the linearization the problem using the Born-approximation~\cite{Kirsch}. In that regime we compute the impulse response for a perturbation of the input  $\eta =  \eta_0 + \delta \eta$, with respect to a reference $\eta_0$. The operator is then accordingly linearized to:
\[
\mathcal{F}^\omega[\eta] 
 = \mathcal{F}^\omega[\eta_0 + \delta \eta]\approx \mathcal{F}^\omega[\eta_0]+F^\omega \delta \eta,
\]
where $F^\omega$ is a linear operator that maps media perturbation to data perturbation. In our case, we assume that $\eta_0=0$, i.e., the background media homogeneous, so $\eta$ itself is the perturbation, and $F$ can be computed explicitly using the free-space Green’s function of the two-dimensional Helmholtz equation. A classical computation omitting the higher order terms in the far-field asymptotics for the Green's function yields:
\begin{equation}\label{eqn:F}
\Lambda^\omega = \mathcal{F}^{\omega} [\eta]\approx F^\omega \eta\,,\quad\text{with}\quad \left (F^\omega \eta \right )(\spmb,\rpmb)=C_\text{nor}\int_{\Omega} e^{-i\omega(\bm{r}-\bm{s})\cdot \bm{y}}\, \eta(\bm{y})\,d\bm{y}\,,
\end{equation}
where $C_\text{nor}= \frac{e^{i\pi/4}}{\sqrt{8\pi\omega}}\omega^2\frac{e^{i\omega R}}{\sqrt{R}}$ is a normalization constant, which we omit from future discussions for the sake of brevity. We stress that according to the formula~\eqref{eqn:F}, the data can be viewed as a first-type Fredholm integration over $\eta$, the perturbation in media. In this perturbed setting, we hope to find $\eta$ that minimize the error between $\Lambda^\omega$ and $F^\omega\eta$:
\begin{equation}
\min_{\eta} \|\Lambda^\omega-F^\omega\eta\|^2, \qquad \text{with} \qquad  
\|\Lambda^\omega-F^\omega\eta\|^2= \int_{\mathbb{S}^1 \times \mathbb{S}^1 }|\Lambda^\omega(\spmb,\rpmb)-\left(F^\omega\eta\right)(\spmb,\rpmb)|^2d \rpmb d\spmb\,. 
\end{equation}
Given that $F^\omega$ is a linear operator on $\eta$, the objective function has a quadratic dependence on the unknown variable $\eta$, and thus, formally, the solution is explicitly given by
\begin{equation}\label{eqn:eta_linear_cont}
\eta
=\left(F^\omega\right)^\dagger\Lambda^\omega\quad\text{with}\quad \left(F^\omega\right)^\dagger=((F^\omega)^\ast F^\omega)^{-1}(F^\omega)^\ast\,,
\end{equation}
where $\left(F^\omega\right)^\ast$ is the adjoint operator of $F^\omega$ that maps data back to the media. This operator is typically called the back-scattering operator.

It can be proved that the operator $(F^\omega)^\ast F^\omega$ is compact thus its numerical inverse will be ill-conditioned. A classical technique to circumvent this issues is to apply the Tikhonov regularization, in which one adds $\epsilon\|\eta\|^2_2$ in the objective function, whose solution is given by
\begin{equation}\label{back-projection}
\eta^\ast=\min_{\eta}\|\Lambda^\omega-F^\omega\eta\|^2+\eps\|\eta\|^2\quad\Rightarrow\quad\eta^\ast = ((F^\omega)^\ast F^\omega+\eps I)^{-1}(F^\omega)^\ast\Lambda^\omega\,.
\end{equation}
Due to this extra regularization term, this strategy is termed the filtered back-projection~\cite{fbp}.

We should stress that one derives this formula~\eqref{back-projection} only in the linearized setting under the assumption that $\eta\sim 0$. It is not valid when $\eta$ is significant, and the formula only serves as the guidance for the actual inversion. In this paper, we are interested in lifting this formula up to the nonlinear setting, and use it as a base to obtain numerical representation of the nonlinear map between data and the media. Mathematically, this is translated to replacing the linear operators~\eqref{back-projection} by nonlinear maps that will be represented by neural networks. In particular, we intend to label $\Lambda^\omega$ and $\eta^\ast$ as the input and output of the neural network, and simulate their relation using a specially designed NN architecture.

Despite~\eqref{back-projection} not being directly applicable, the formula nevertheless provides insightful guidelines in designing neural architectures. It suggests that the relation between the data and the media is composed of two operations : the inner-layer operator is the adjoint operator $(F^\omega)^*$ (or the back-scattering operator), and the outer-layer operator, $(F^\omega)^\ast F^\omega+\eps I$ (or filtering operator). When one designs neural networks to represent the inverse operator \eqref{eq:inv_formulation}, both operators need to be integrated. Each enjoy some unique features, which should be respected and therefore encoded in the neural architecture. In particular, the filtering operator is a pseudo-differential operator of convolution type, therefore \textit{translational equivariant}. As a consequence, the neural ansatz is also of convolution-type and should respect this translational equivariance property. Indeed, we choose to represent it using a few two-dimensional convolutional layers as its discretization, as was studied in~\cite{MLZ,FanYing:scattering}.

The presentation of the back-scattering operator $(F^\omega)^*$ in the nonlinear setting is more involved. One can prove that this operator satisfies the \textit{rotational equivariance}. This property should be respected by its neural representation. Furthermore, this operator satisfies the \textit{complementary low-rank property}, thus it is expected to be well approximated by the butterfly factorization. Accordingly, the neural representation could incorporate the butterfly structure as well. In Section~\ref{sec:propertyequivariance} and~\ref{sec:propertybutterfly} respectively, we will study in depth these two properties, and investigate how to encode these features in the NN architecture to achieve both inversion accuracy and efficiency through enforcing this particular type of compressibility.

\subsection{Discretization}\label{sec:setup}

We translate the discussion from the previous sections to the discrete setting. To streamline the notation, quantities in {calligraphic} fonts, such as $\mathcal{F}^\omega$, are used to denote nonlinear maps, while the ones in regular fonts, such as $F^\omega$ and $\Lambda^\omega$, are used to denote the linearized version. The quantities written in serif font, such as $\sfF^\omega$ and $\sf\Lambda^\omega$, are used to present the discretized version of the associated linear operators.

Our method relies on particular choices of the discretization, which we described in what follows. Throughout the paper, we parametrize the incoming direction and the sampling point, $\bm{s},\bm{r}\in\mathbb{S}^1$, by their associated angles
\[
\bm{s} = (\cos(s), \sin(s))\text{\quad and\quad}\bm{r} = (\cos(r), \sin(r))\,,
\]
and correspondingly, we denote the continuous normalized far-field pattern as 
\begin{equation}\label{eqn:normalized}
    \Lambda_\nor^\omega(r,s) := u^\sca(R\rpmb;\spmb)/C_\nor\,,
\end{equation}

and the normalized operator $F^\omega$ by $F^\omega_\nor$. Upon the normalization, the equation~\ref{eqn:F} becomes 
\begin{equation}
\Lambda_\nor^\omega= F_\nor^\omega\eta\,,\quad\text{with}\quad
\Lambda_\nor^\omega(r,s)
=\int e^{-i\omega(\bm{r}-\bm{s})\cdot \bm{y}}\eta(\bm{y})\,d\bm{y}\,.
\end{equation}
which leads to the filtered back-projection
\[\eta^\ast = ((F_\nor^\omega)^\ast F_\nor^\omega+\eps I)^{-1}(F_\nor^\omega)^\ast \Lambda_\nor^\omega\,.\]

For the conciseness of the notation, we drop the subscript $\cdot_\nor$. Numerically, the directions of sources and detectors are represented by the same uniform grid in $\mathbb{S}^1$ with $n_{\text{sc}}$ grid points given by
\[
s_j, r_j = \frac{2\pi j}{n_{\text{sc}}},\ j=0,\dots,n_{\text{sc}}-1\,.
\]
Using this setting, the discrete data $\sf\Lambda^\omega$ takes its values on the tensor product of both grids with complex values, which are decomposed in their real and imaginary parts
\begin{equation}\label{eqn:def_dis_d}
\sf\Lambda^\omega=\Lambda^\omega_R+i\Lambda^\omega_I\in \C^{n_\sca\times n_\sca}\,. 
\end{equation}

{For the discretization of $\Omega$, we use two discretizations, each of them anchored on specific properties of the operators, which will be justified in Section~\ref{sec:propertyequivariance}.}
On the one hand, as it will be shown in the sequel, the back-scattering operator $(F^\omega)^\ast$ enjoys the rotational equivariance property (Section~\ref{sec:propertyequivariance}), 
and it can naturally be represented on the polar coordinates. On the other hand, the filtering operator $((F^\omega)^\ast F^\omega+\eps I)^{-1}$ is translational equivariant and it should be represented on a Cartesian grid.

Therefore, we set the physical domain to be $\Omega=[-0.5,0.5]^2$. We use a Cartesian mesh of $n_\eta\times n_\eta$ grids. We also discretize it using  polar coordinates on a slightly larger domain. In particular, we note $\Omega\subset B_{1/2}(0)$, and define
\begin{equation}\label{eqn:theta}
    (\theta_j,\rho_i) = \bigg(\frac{2\pi j}{n_{\theta}}\,,\frac{i}{2n_{\rho}}\bigg)\,,\quad\text{with}\quad j=0,\dots,n_{\theta}-1\,, i=0,\dots,n_{\rho}-1\,.
\end{equation}

Thus, $\eta(\pmb{x})$ is represented as a tensor: $\eta\in\R^{n_\eta\times n_\eta}$ with its values being $\eta(\pmb{x})$ evaluated on the Cartesian mesh. We also assume $n_\theta=n_\sca$ so that $s_j,r_j$ and $\theta_j$ are sampled on the same mesh.


As a consequence, in the algorithmic pipeline we need to express the intermediate representations in both polar and Cartesian meshes on the physical domain, which in return requires a change of coordinates of the data between the application of the two operators. We denote the discrete intermediate field obtained upon the application of the discretized back-scattering operator $(\sfF^\omega)^\ast$ by $\alpha^\omega$:
\begin{equation}\label{eqn:def_alpha}
\sf\alpha^\omega\defeq(F^\omega)^\ast \Lambda^\omega\,,
\end{equation}
It lives in the range of $(\sfF^\omega)^\ast$ and is represented in polar coordinates. Thus we sample $\sf\alpha^\omega$ on a grid using the polar coordinates, meaning $\alpha^\omega$ is presented as a vector of size $n_\theta \cdot n_\rho$. Then we use quadratic interpolation, to interpolate it to a Cartesian mesh before being fed to the filtering operator
\begin{equation}\label{eqn:def_eta_alpha_relation}
\sf\eta = ((F^\omega)^\ast F^\omega+\eps I)^{-1}\alpha^\omega\,.
\end{equation}

\subsection{Computation of \texorpdfstring{$(F^\omega)^*$}{F*}}\label{sec:F_ast}
We claim that this back-scattering operator enjoys both compressibility through the butterfly structure and the rotational equivariance property. To do so, we compute an explicit formulation for this operator in what follows. Noting that using the standard inner product, we can flip $(F^\omega)^\ast$ to its dual space:
\begin{align*}
    \inner{\eta,\alpha^\omega}_{\Omega}=\inner{\eta,(F^\omega)^*\Lambda^\omega}_{\Sbb^1\times\Sbb^1}\,.
\end{align*}
With straightforward calculation, recalling~\eqref{eqn:F}, the right hand side becomes:
\begin{equation}
\begin{aligned}
    \inner{\eta,\alpha^\omega}_{\Omega}=\inner{\eta,(F^\omega)^*\Lambda^\omega}_{\Omega}
    &=\int\barr{\Lambda^\omega(r,s)}\int e^{-i\omega(\bm{r}-\bm{s})\cdot \bm{y}}\eta(\bm{y})\,d\bm{y}\,ds\,dr\,,\\
&=\int \eta(\bm{y})\barr{\int e^{i\omega(\bm{r}-\bm{s})\cdot \bm{y}} \Lambda^\omega(r,s)\,ds\,dr}\,d\bm{y}\,,
\end{aligned}
\end{equation}

suggesting
\begin{equation}\label{eqn:F_ast}
(F^\omega)^* \Lambda^\omega(\bm{y})=\int e^{i\omega(\bm{r}-\bm{s})\cdot \bm{y}}\Lambda^\omega(r,s)\,ds\,dr\,.
\end{equation}

To view it in polar coordinates, we denote that $\bm{y} = (\rho\cos\theta,\rho\sin\theta)$ and have 
\[
    i\omega(\bm{r}-\bm{s})\cdot \bm{y} = i\omega\rho(\cos(r-\theta)-\cos(s-\theta))\,.
\]
Injecting it back to the formula above, the calculation becomes\footnote{$\Lambda^{\omega}$ is extended periodically outside the domain $(r,s)\in[0,2\pi]^2$.}
\begin{equation}\label{FIO}
\begin{aligned}
\alpha^\omega(\theta,\rho)=((F^\omega)^* \Lambda^\omega)(\theta,\rho) &= \iint_{[0,2\pi]^2} e^{i\omega\rho\cos(r-\theta)}e^{-i\omega\rho\cos(s-\theta)}\Lambda^\omega(r,s)\,ds\,dr\,,\\
&= \iint_{[0,2\pi]^2} e^{i\omega\rho\cos(r)}e^{-i\omega\rho\cos(s)}\Lambda^\omega(r+\theta,s+\theta)\,ds\,dr\,,\\
&= \int_{[0,2\pi]} e^{i\omega\rho\cos(r)}\bigg(\int_{[0,2\pi]}e^{-i\omega\rho\cos(s)}\Lambda^\omega(r+\theta,s+\theta)\,ds\bigg)\,dr\,.
\end{aligned}
\end{equation}

The formula suggests that the adjoint operator $\left(F^\omega\right)^\ast$ can be decomposed in two-integral operators, both of which involves the integration kernel 
\begin{equation}\label{eqn:kernal}
    K^\omega(\rho,t) \defeq e^{-i\omega\rho\cos(t)}\,.
\end{equation}

\subsection{Equivariance}\label{sec:propertyequivariance}
Both the back-scattering operator and the filtering operator in~\eqref{back-projection} enjoy certain equivariances and the design of NN should respect these features. We summarize the two equivalance that will be incorporated in our studies.
\paragraph{Rotational equivariance} One key feature that $(F^\omega)^\ast$ enjoys is the rotational equivariance: It means that by rotating the data by a certain angle, one is expected to obtain an output that is rotated by the same angle. Recalling the definition of the data~\eqref{eqn:def_alpha}, and defining the following operator: 
\begin{equation}
    \mathcal{R}_{\tilde\theta} [\Lambda^\omega](r,s) = \Lambda^\omega(r-\tilde\theta,s-\tilde\theta)\,,
\end{equation}
which rotates by $\tilde\theta$ the data. Then we expect 
\begin{equation}\label{eqn:def_equi}
    (F^\omega)^\ast \mathcal{R}_{\tilde\theta} [\Lambda^\omega] (\theta, \rho) =  \mathcal{R}_{\tilde\theta}[(F^\omega)^\ast \Lambda^\omega]  (\theta, \rho) = [(F^\omega)^\ast \Lambda^\omega](\theta - \tilde\theta, \rho)\,,
\end{equation}
in which the operator only acts on the angular inputs. This property can be easily checked by examining the explicit expression~\eqref{FIO}. Considering that neither integral kernel has $\theta$ dependence and by trivially rearranging terms:
\[
\begin{aligned}
(F^\omega)^\ast \mathcal{R}_{\tilde\theta} [\Lambda^\omega] (\theta, \rho) & =  \int_{[0,2\pi]} e^{i\omega\rho\cos(r)}\bigg(\int_{[0,2\pi]}e^{-i\omega\rho\cos(s)}\Lambda^\omega((r - \tilde\theta)+\theta,(s-\tilde\theta)+\theta)\,ds\bigg)\,dr, \\
& =  \int_{[0,2\pi]} e^{i\omega\rho\cos(r)}\bigg(\int_{[0,2\pi]}e^{-i\omega\rho\cos(s)}\Lambda^\omega(r+(\theta-\tilde\theta),s+(\theta- \tilde\theta))\,ds\bigg)\,dr, \\
& = (F^\omega)^\ast [\Lambda^\omega] (\theta - \tilde \theta, \rho) = \mathcal{R}_{\tilde\theta}[(F^\omega)^\ast \Lambda^\omega]  (\theta, \rho)\,,
\end{aligned}
\]
proving~\eqref{eqn:def_equi}.

\paragraph{Translational equivariance} 

Using the explicit formulas of $F^\omega$ in~\eqref{eqn:F} and $(F^\omega)^\ast$ in~\eqref{eqn:F_ast}, we can also verify the translational equivariance of $((F^\omega)^\ast F^\omega+\eps I)^{-1}$. Using the Fubini's theorem, we have
\begin{equation}
\begin{aligned}
    (F^\omega)^\ast F^\omega(\eta)(\bm{y}) &= \int_{[0,2\pi]^2} e^{i\omega(\bm{r}-\bm{s})\cdot \bm{y}} \left (\int_{\Omega} e^{-i\omega(\bm{r}-\bm{s})\cdot \bm{x}}\, \eta(\bm{x})\,d\bm{x} \right ) \, \,ds\,dr\,,\\
    &= \int_{[0,2\pi]^2}\int_{\mathbb{R}^2} e^{i\omega(\bm{r}-\bm{s})\cdot \bm{y}} e^{-i\omega(\bm{r}-\bm{s})\cdot \bm{x}}\, \eta(\bm{x})\,d\bm{x}\, \,ds\,dr\,,\\
    &= \int_{\mathbb{R}^2} \left (\int_{[0,2\pi]^2}e^{i\omega(\bm{r}-\bm{s})\cdot (\bm{y}-\bm{x})} \,ds\,dr\,\right ) \, \eta(\bm{x})\,d\bm{x}\,,\\
    & = \int p(\bm{y}-\bm{x})\eta(\bm{x})d\bm{x} = p\ast \eta(\bm{y}).\\
\end{aligned}
\end{equation}
where the convolution kernel is
\begin{equation}
    p(\bm{x}) := \int_{[0,2\pi]^2}e^{i\omega(\bm{r}-\bm{s})\cdot \bm{x}} \,ds\,dr\,,
\end{equation}
and the function $\eta$ extended by zero on $\Omega^c$. This convolution feature justifies that $(F^\omega)^\ast F^\omega$ is translational equivariant. The same property also holds when we add it with an identity, and take its inversion.

\subsection{Compression of the kernel \texorpdfstring{$\sfK^\omega$}{K}}\label{sec:propertybutterfly}
The back-scattering operator~\eqref{FIO} presents two levels of integration against the oscillatory integral kernel $K^\omega(\rho,t)=e^{-i\omega\rho\cos(t)}$ (or its complex conjugate). A brute-force integration of such an oscillatory integrand typically calls for very fine discretization and, therefore, high cost. However, following~\cite{MLZ}, one can argue that its associated discretized form $\sfK^\omega$ can be compressed through the butterfly factorization~\cite{BF} by leveraging the complimentary low-rank property~\cite{BF}. In this section we provide an overview of (and the intuition behind) the butterfly factorization algorithm.

\begin{figure}[h]
\centering
\includegraphics[scale = 0.3]{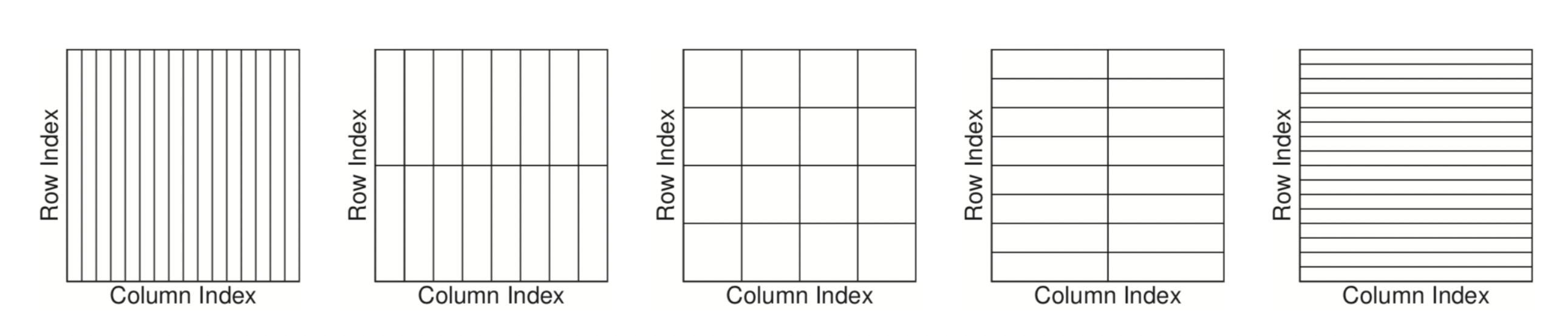}
\caption{Sketch of a family of partitions of a matrix exhibiting the complementary low-rank property. Each sub-matrix induced by the different partition has the same numerical rank. }\label{fig:bf}
\end{figure}

For one-dimensional problems the complementary low-rank property can be stated as follows: any block of the matrix with constant area, i.e., the multiplication of its heights by its width, has its numerical rank upper bounded by a constant. For example, in the family of partitions depicted in Figure~\ref{fig:bf} each block has the same area, and therefore in order to satisfy the complimentary low-rank condition, each block should have the same numerical rank.

We claim that the integrating kernel $K^\omega(\rho,t)=e^{-i\omega\rho\cos(t)}$ satisfies complementary low-rank property. Indeed, the claim can be broadened to any function in the form of $e^{i\omega\phi(\rho,t)}$. To see it, we perform Taylor expansion in the neighborhood of $(\rho_0,t_0)$. Let $\abs{\rho-\rho_0}<d_\rho$  and $\abs{t-t_0}<d_t$, we can approximate $\phi(\rho,t)$ up to the second order expansion as
\begin{equation}
\begin{aligned}
\phi(\rho,t) &= \phi(\rho_0,t_0) + \del_\rho \phi(\rho_0,t_0) \cdot (\rho-\rho_0) + \del_t \phi(\rho_0,t_0) \cdot (t-t_0)\\
&+ (\rho-\rho_0)^T\cdot\del^2_\rho \phi(\rho_0,t_0) \cdot (\rho-\rho_0) + (t-t_0)^T\cdot\del^2_t \phi(\rho_0,t_0)\cdot (t-t_0)+\mathcal{O}(d_\rho d_t)\,.
\end{aligned}
\end{equation}

Since the first five terms are separable and the reminder term is bounded by $\mathcal{O}(d_\rho d_t)$, we have
\begin{equation}\label{eqn:butterfly_expansion}
e^{i\omega\phi(\rho,t)} = e^{i\omega\psi(\rho)}e^{i\omega\xi(t)}e^{i\omega\mathcal{O}(d_\rho d_t)} = e^{i\omega\psi(\rho)}e^{i\omega\xi(t)}(1+\mathcal{O}(\omega d_\rho d_t))\,,
\end{equation}
where $\psi(\rho)= \phi(\rho_0,t_0) + \del_\rho \phi(\rho_0,t_0) \cdot (\rho-\rho_0)+ (\rho-\rho_0)^T\cdot\del^2_\rho \phi(\rho_0,t_0) \cdot (\rho-\rho_0)$ and $\xi(t)=\del_t \phi(\rho_0,t_0) \cdot (t-t_0)+(t-t_0)^T\cdot\del^2_t \phi(\rho_0,t_0)\cdot (t-t_0)$. Therefore, when $d_\rho d_t<\omega^{-1}$, $e^{i\omega\phi(\rho,t)}$ can be locally approximated by separable functions. Further details can be found in~\cite{MLZ}.

We note that the error term from a separable form is a product of two terms, the size of interval in $\rho$ and the size of interval in $t$. If the product of these two distances is small, the integral kernel can be approximated by a separable function, hence becoming low-rank. This 
 observation is similar to the requirement of the complementary low-rank property where the size of column indexing and that of row indexing in a given block of the matrix need to complement each other, so the total number of entries (or its ares) in this block is controlled.

A matrix satisfying the complementary low-rank property can be expanded on both sides using a factorization in highly-structured sparse matrices, which can be written using 
smaller-sized matrix-matrix products with operations on rows/columns respectively, as illustrated in the Figure~\ref{fig:bf}. This multiplication strategy is often termed \textit{butterflying} a matrix. It allows a lower complexity matrix-vector product, reducing the cost of from $N^2$ to $\mathcal{O}(N\log{N})$~\cite{BF}. In our case, it is $\sfK$, a matrix of size $N\times N$ that is complementary low-rank and thus is butterfliable~\cite{BF}. To be more specific, the butterfly algorithm should find an expansion of the matrix by a product of $L+3$ sparse matrices, each of which having only $\mathcal{O}(N)$ entries:

\begin{equation}\label{eqn:butterfly}
\sfK\approx \sfU^L\sfG^{L-1}\cdots \sfG^{L/2}\sfM^{L/2}(\sfH^{L/2})^*\cdots(\sfH^{L-1})^*(\sfV^L)^*\, .
\end{equation}

This formula approximately expand $\sfK$ into the multiplication of a series of matrix products. Denoting $L$ the number of ``levels'' in the factorization, $\sfU^L$ and $\sfV^L$ will be block diagonal matrices, and $\sfM^{L/2}$ is a weighted permutation matrix, which is usually called a switch matrix. The structures of the factors in the butterfly factorization provide acute intuition into their interpretations. For example, when a vector is right multiplied by the matrix, $\sfU^L$ extracts a local representation of the vector, and then each $\sfG^l$ compresses two adjacent local representations. Upon the application of the switch matrix $\sfM^{L/2}$ that redistributes the representations from the previous step by permuting the vector, each $\sfH^l$ decompress the representation by splitting it into two, which increases the resolution of the representation. Finally, $\sfV^L$ converts the local representations to sampling points. 

Algorithmically, this decomposition is achieved in two stages. In the first, one would perform a singular value decomposition (SVD) for submatrices at level $L/2$. This step is composed of SVDs for $N$ submatrices, with each submatrix of size $m\times m$, where $m=\sqrt{N}$. For the moment, we unify the rank for each SVD decomposition to be $r$. Denoting $\sfU^{L/2}_{i,j}\in\mathbb{R}^{m\times r}$, $\sfS^{L/2}_{i,j}\in\mathbb{R}^{r\times r}$, and $\sfV^{L/2}_{i,j}\in\mathbb{R}^{m\times r}$ the left singular vectors, singular values and right singular vectors, we stack them up as
\begin{equation}
\begin{aligned}
    \sfK&\approx\sfU^{L/2}\sfM^{L/2}(\sfV^{L/2})^\ast\\
    =&\begin{pmatrix}
\sfU_{0,0}^{L/2}\sfS_{0,0}^{L/2}(\sfV^{L/2}_{0,0})^\ast & \sfU_{0,1}^{L/2}\sfS_{0,1}^{L/2}(\sfV^{L/2}_{1,0})^\ast & \dots & \sfU_{0,m-1}^{L/2}\sfS_{0,m-1}^{L/2}(\sfV^{L/2}_{m-1,0})^\ast\\
\sfU_{1,0}^{L/2}\sfS_{1,0}^{L/2}(\sfV^{L/2}_{0,1})^\ast & \sfU_{1,1}^{L/2}\sfS_{1,1}^{L/2}(\sfV^{L/2}_{1,1})^\ast &  & \sfU_{1,m-1}^{L/2}\sfS_{1,m-1}^{L/2}(\sfV^{L/2}_{m-1,1})^\ast\\
\vdots& &\ddots& \\
\sfU_{m-1,0}^{L/2}\sfS_{m-1,0}^{L/2}(\sfV^{L/2}_{0,m-1})^\ast & \sfU_{m-1,1}^{L/2}\sfS_{m-1,1}^{L/2}(\sfV^{L/2}_{1,m-1})^\ast & \dots & \sfU_{m-1,m-1}^{L/2}\sfS_{m-1,m-1}^{L/2}(\sfV^{L/2}_{m-1,m-1})^\ast
\end{pmatrix}\,.
\end{aligned}
\end{equation}
The approximation sign takes into account that SVD cuts off small singular values, and
\begin{equation}\label{eqn:stack_U}
\sfU^{L/2} = \begin{pmatrix}
    \sfU_0^{L/2} & & & \\
& \sfU_1^{L/2} & &\\
& & \ddots & \\
& & &  \sfU_{m-1}^{L/2}
\end{pmatrix}\quad\text{and}\quad
(\sfV^{L/2})^\ast = \begin{pmatrix}
(\sfV_0^{L/2})^\ast & & & \\
&(\sfV_1^{L/2})^\ast & &\\
& & \ddots & \\
& & &(\sfV_{m-1}^{L/2})^\ast
\end{pmatrix}\,.
\end{equation}
with
\begin{equation}\label{eqn:stack_U_i}
\sfU_i^{L/2} = \begin{pmatrix}
\sfU_{i,0}^{L/2}& \sfU_{i,1}^{L/2}&\dots&\sfU_{i,m-1}^{L/2}
\end{pmatrix}\quad\text{and}\quad\sfV_i^{L/2} = \begin{pmatrix}
\sfV_{i,0}^{L/2}& \sfV_{i,1}^{L/2}&\dots&\sfV_{i,m-1}^{L/2}
\end{pmatrix}\,.
\end{equation}

Note that the diagonal blocks of $\sfU$, denoted as $\sfU_i$ stores left (column) singular vectors for the $i$-th row block, and the total size of $\sfU^{L/2}_i$ is $m\times mr$. Similarly, diagonal blocks of $\sfV$ stores right (column) singular vectors of size $m$. In the second stage, one expands these vectors in $\sfU$ and $\sfV$ recursively, by partitioning them in half, while expanding them in the opposite direction to connect parallel blocks. In particular, for each diagonal block $\sfU_i$, we split it in top and bottom halves
\[
\sfU^{L/2}_i=\begin{pmatrix}
\sfU^{L/2,t}_{i}\\
\sfU^{L/2,b}_{i}
\end{pmatrix} = \begin{pmatrix}
\sfU^{L/2,t}_{i,0} & \sfU^{L/2,t}_{i,1} & \cdots  & \sfU^{L/2,t}_{i,m-1}\\
\sfU^{L/2,b}_{i,0} & \sfU^{L/2,b}_{i,1} & \cdots  & \sfU^{L/2,b}_{i,m-1}\\
\end{pmatrix} \,.
\]
Noting that $\left(\sfU^{L/2,t}_{i,2j},\sfU^{L/2,t}_{i,2j+1}\right)$ form the column space of $\sfK^{L/2+1}_{2i,j}$, the $(2i,j)$-block of $\sfK$ decomposed at $(L/2+1)$-th level for row, and $(L/2-1)$-th level for column. Similarly $\left(\sfU^{L/2,b}_{i,2j},\sfU^{L/2,b}_{i,2j+1}\right)$ form the column space of $\sfK^{L/2+1}_{2i+1,j}$, there are translation matrix $\sfG^{L/2}_{2i,j}$ and $\sfG^{L/2}_{2i+1,j}$ so that:
\[
\left(\sfU^{L/2,t}_{i,2j},\sfU^{L/2,t}_{i,2j+1}\right)=\sfU^{L/2+1}_{2i,j}\sfG^{L/2}_{2i,j}\,,\quad
\left(\sfU^{L/2,b}_{i,2j},\sfU^{L/2,b}_{i,2j+1}\right)=\sfU^{L/2+1}_{2i+1,j}\sfG^{L/2}_{2i+1,j}\,.
\]
Using the same definition as in~\eqref{eqn:stack_U} and~\eqref{eqn:stack_U_i} to stack up the matrices, we accordingly define $\sfU^{L/2+1}_i$ and $\sfU^{L/2+1}$. These definitions also allow us to write the transform in a concise form with properly stacked up $\sfG^{L/2}$:
\[
\sfU^{L/2} = \sfU^{L/2+1}\sfG^{L/2}\,.
\]

Viewing these transformation, it is clear that the row size of $\sfU^{L/2+1}_{2i,j}$ is that of $\sfU^{L/2,t}_{i,2j}$, which is half of that of $\sfU^{L/2}_{i,2j}$. Similarly, $\sfU^{L/2+1}_{2i,j}$ has the same rank $r$, and thus $\sfU^{L/2+1}_{2i,j}$ is of size $\frac{m}{2}\times\frac{mr}{2}$. As a summary, when one changes from $L$-th iteration to $L+1$, the number of blocks double in both row and column, while the rank of each sub-block is kept as $r$.

\begin{figure}[h]
\centering
\includegraphics[scale = 0.3]{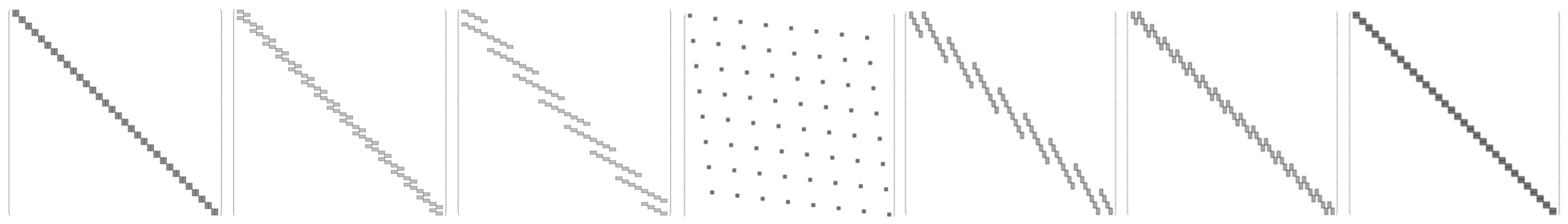}
\caption{An illustration of the matrix factors in the butterfly factorization. In the illustration, $L=6$, $r=1$, and $N = 64$. }\label{fig:BF}
\end{figure}

Perform this iteration recursively by $L/2$ steps, we finally arrive at the formulation of~\eqref{eqn:butterfly}. At the final stage of the iteration, $\sfU^L$ as defined in the same way as in~\eqref{eqn:stack_U} is composed with $N$ blocks along its diagonal, with each block being of $1\times r$. As an illustration, a rank $1$ approximation using the butterfly algorithm can shown visualized in Fig.~\ref{fig:BF}.

\section{Neural Network Architecture}
We design a neural architecture to conduct the inverse scattering in this section.

As discussed in the introduction, the implementation of $(F^\omega)^\ast$ and $((F^\omega)^\ast F^\omega+\eps I)^{-1}$ should honor the operators' own properties. In particular, we used convoluational NN to lift the filtering operator so it satisfies the translational equivariance. As argued before the back-scattering operator has rotational equivariance, and it is possible to compress it via a butterfly structure. We will discuss these properties in this section, and present how they get integrated in the design of NN in Section~\ref{sec:equivariance} and~\ref{sec:butterfly} respectively.

We write out architecture as 
\begin{equation}
\eta = \Phi_{\Theta} (\{\Lambda^{\omega}\}_{\omega \in \bar{\Omega}})\,,
\end{equation}
where the input  $\{\Lambda^{\omega}\}_{\omega \in \bar{\Omega}}$ is a collection of far-field patterns as defined in \eqref{eq:far_field_pattern}, indexed by frequencies in $\bar{\Omega}$. $\Phi_\Theta$ is the function generated by the neural network with neurons weighted using parameters $\Theta$. 

If we have only one frequency, i.e., $\bar{\Omega} = \{\omega\}$, then our architecture follows closely the filtered back-projection with a change of variables as described in Alg.~\ref{alg:filtered_back_projection}. In practice we modify the application of $F^{\omega}$ by a neural network that mimics the butterfly structure, and we replace the filtering by several layers of convolutional networks.

\begin{algorithm}[H]
\begin{algorithmic}[1]
\caption{Filtered Back-Projection.} \label{alg:filtered_back_projection}
\Statex Input: {$\sf\Lambda^\omega \in \mathbb{C}^{n_{\sca} \times n_{\sca}} $}
\Statex Output: {$\sf\eta = (( F^\omega)^* F^\omega + \epsilon I)^{-1}(F^\omega)^*$}{$\Lambda^\omega \in \mathbb{R}^{n_{\eta} \times n_{\rho}}$}
    \item Preparation of the normalized data~\eqref{eqn:normalized};
    \item Computation of back-projection  $\sf\alpha^\omega(\theta,\rho) = (F^\omega)^\ast \Lambda^\omega$ using \eqref{eqn:def_alpha};
    \item Change of variables of $\alpha^\omega$ to Cartesian coordinates $\alpha^\omega(\theta,\rho) \rightarrow \alpha^\omega(x,y)$ through polynomial interpolation;
    \item Computation of the filtering operation  $\sf\eta = ((F^\omega)^\ast F^\omega+\eps I)^{-1} \alpha^\omega$ using~\eqref{eqn:def_eta_alpha_relation}.
\State \textbf{return} $\eta $
\end{algorithmic}
\end{algorithm}

\begin{figure}[h]
    \centering
    \includegraphics[width=1\textwidth]{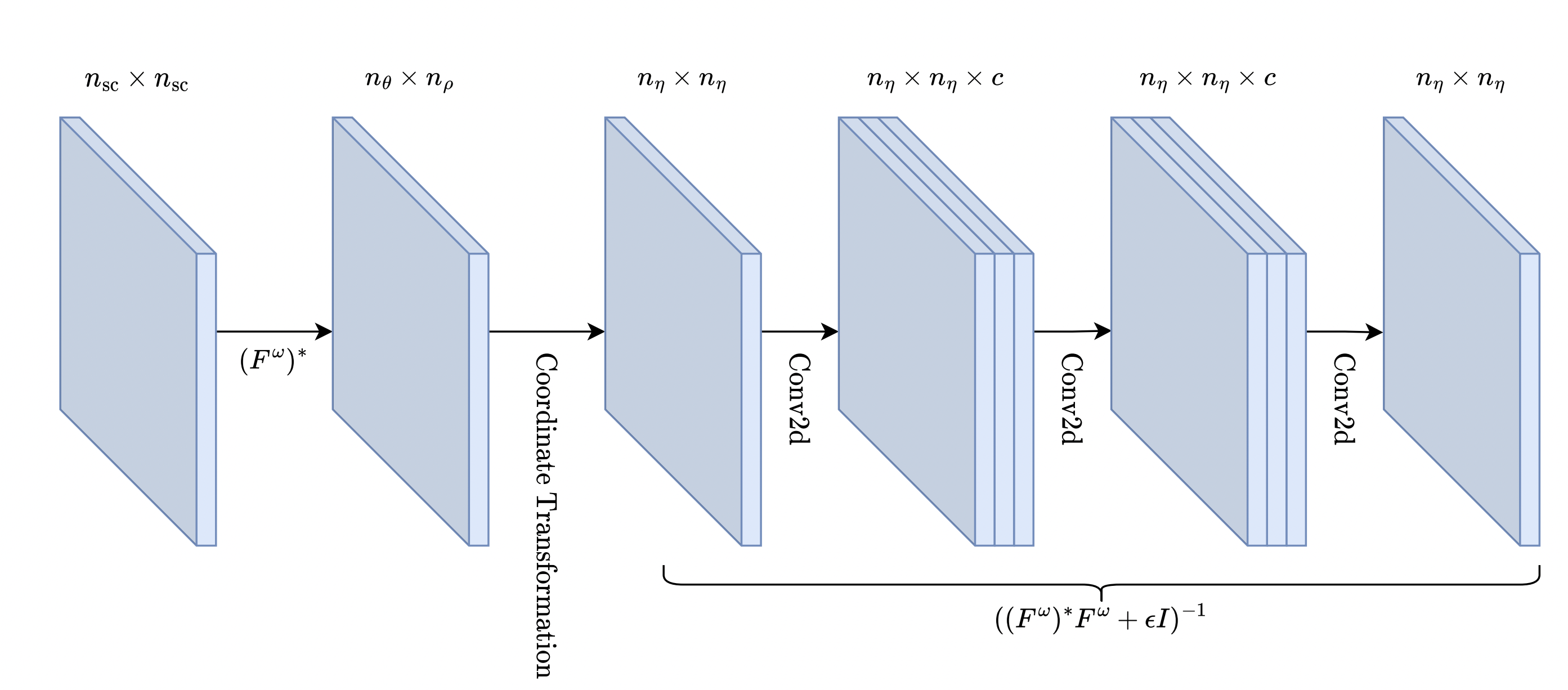}
    \caption{An illustration of the architecture of the model. The shape of each tensor at each step is labelled on the top. Convolutional NN is used to lift the filtering operator.}
    \label{fig:procedure}
\end{figure}

{For the multiple frequencies case, we consider, as a motivation, the loss function
\begin{equation}
    \eta^{*} = \text{argmin} \sum_{\omega \in \bar{\Omega}} \| \Lambda^{\omega} - F^{\omega}[\eta]  \|^2 + \epsilon \|\eta \|^2,
\end{equation}
and by computing the first variation we have that 
\begin{equation}
    \eta^{*} = \sum_{\omega \in \bar{\Omega}}  (( F^\omega)^* F^\omega + \epsilon I)^{-1}(F^\omega)^*\Lambda^\omega.
\end{equation}}



This can be easily implemented by applying Alg.~\ref{alg:filtered_back_projection} to each component and then average them at the end. However, as it was shown in \cite{MLZ} this approach does not mix the different scales efficiently. One alternative, would be to use the multilevel structure of the butterfly factorization of $K^{\omega}$ and merge the different scales hierarchically when applying the neural equivalent of $F^{\omega}$. The implementation of this pipelines is however cumbersome and prone to errors. Instead we considered an intermediate architecture: we apply the \textit{back-projection} independently for each frequency, we concatenate the resulting $\alpha^{\omega}$'s along a new dimension, equivalent to a channel dimension for CNN. Finally, we use a \textit{shared} filter acting in all frequencies together, so the mixing and weighting for data at different frequencies can be done progressively.

In a nutshell, for the general case, our neural network initially incorporates data of different frequencies when computing the back-projections. The $\alpha^\omega$ for different frequencies $\omega$ are concatenated along a channel dimension, and they are merged when using a global filter operation. The method is summarized in Alg.~\ref{alg:wide-band-equivariant-network}. 


\begin{algorithm}[H]
\begin{algorithmic}[1]
\caption{Wide-Band Equivariant Network.} \label{alg:wide-band-equivariant-network}
\Statex Input: {$\{\sf\Lambda^\omega \}_{\omega \in \bar{\Omega}}$, where $ \sf\Lambda^\omega \in \mathbb{C}^{n_{\sca} \times n_{\sca}},\,\, \text{for } \omega \in \bar{\Omega} $}
\Statex Output: {$\sf\eta \in \mathbb{R}^{n_{\eta} \times n_{\rho}}$}
{
\For{$\omega \in \bar{\Omega}$}           
         \State Preparation of the normalized data~\eqref{eqn:normalized};
    \State Computation of back-projection  $\sf\alpha^\omega(\theta,\rho) = (F^\omega)^\ast \Lambda^\omega$ using \eqref{eqn:def_alpha};
    \State Change of variables of $\alpha^\omega$ into Cartesian coordinates $\alpha^\omega(\theta,\rho) \rightarrow \alpha^\omega(x,y)$ through polynomial interpolation;
\EndFor
   }
\State {Concatenation of $\{\alpha^\omega\}_{\omega \in \bar{\Omega}}$ into a 3d tensor $\mathbf{a}$ such that $\mathbf{a} [\,:,:,\omega] = \alpha^\omega$ .}
    \item Computation of the filtering operation  $\sf\eta = ((F^\omega)^\ast F^\omega+\eps I)^{-1} \mathbf{a}$ using~\eqref{eqn:def_eta_alpha_relation}.
\State \textbf{return} $\eta $
\end{algorithmic}
\end{algorithm}

{In computation, the size and the values in the set of frequencies $\bar{\Omega}$ depends on the target resolution. We follow a dyadic partition for the frequencies. For instance, we will use data corresponding to source frequencies 2.5, 5.0, and 10.0 Hz for data of dimension $n_\sca=80$. Hence, the number of frequencies scale logaritmically with respect to $n_\sca$.}

\subsection{Application of equivariance}\label{sec:equivariance}
Recall the formula~\eqref{FIO}, the backscattered data $\alpha^\omega$ comes from two layers of integral operators, both of which involves the integration kernel $K^\omega(\rho,t) \defeq e^{-i\omega\rho\cos(t)}$. Using the discretization introduced in Section~\ref{sec:setup}, this function is presented using the matrix form
as
\[
\mathsf{K}^\omega \in \R^{n_\sca\times n_\rho}\,,\quad\text{with}\quad \mathsf{K}^\omega_{mn} = e^{-i\omega\rho_n\cos(t_m)} = K^\omega(\rho_n,t_m)\,.
\]
Recalling the discrete form of data~\eqref{eqn:def_dis_d}, to fully implement the back-scattering operator~\eqref{FIO}, one can rewrite~\eqref{FIO}, for any $\theta_j$~\ref{eqn:theta}, as:
\begin{equation}
\label{eqn:two_forms}
\alpha^\omega(\theta_j,\cdot)=\sf((F^\omega)^*\Lambda^\omega)(\theta_j,\cdot) = \underbrace{\operatorname{ones}(1,n_\sca)\cdot
[\barr{\mathsf{K}^\omega} \odot (\sf\Lambda^\omega_{\theta_j}\cdot \mathsf{K}^\omega )]}_{\text{{Implementation} I}} = \underbrace{\operatorname{diag}
[(\mathsf{K}^\omega)^\ast \cdot \sf\Lambda^\omega_{\theta_j}\cdot \mathsf{K}^\omega]}_{\text{{Implementation} II}}\,.
\end{equation}

Here $\cdot$ denotes the matrix multiplication and $\odot$ denotes the element-wise Hadamard multiplication. $\cdot^\ast$ is the notation for  conjugate transpose. $\sf\Lambda^\omega_\theta$ is the discrete shifting of $\Lambda^\omega(r+\theta,s+\theta)$. Considering $\{\theta_j\}=\{r_j\}$, for the matrix, we directly shift all rows/columns by $j$:
\begin{equation}\label{eqn:data_shift}
\sf\Lambda^\omega_{\theta_j}(m,n)=\sf\Lambda^\omega(m+j,n+j)\,.
\end{equation}
It should be noted that the two {implementations} are mathematically equivalent, although they are numerically different. The compression of $\sfK^\omega$ drives the choice of the implementation. 
In particular, Implementation I replaces the first matrix-matrix-product using a Hadamard entry-wise product so it reduces the  computational cost,
and it is therefore preferred when there is no compression. On the other hand, this implementation is unfriendly to encode the butterfly structure. The butterfly structure requires the to-be-examined matrix to be expanded from both sides into many smaller-sized matrix products, and this is incompatible to the Hadamard product. As a consequence, when the  butterfly-structure is used for compression, we turn to Implementation II.

Both implementations impose the integral kernel $\sfK^\omega$ to be shared across different $\theta_j$. 
This observation lies at the core of our algorithm.
Its advantages are twofold: it automatically enforces the equivariance property encoded in~\eqref{eqn:def_equi}, and it significantly reduces the computational cost. 

Indeed, with this approach, the to-be-trained parameters are all encoded in $\sfK^\omega$, a small matrix of size ($n_\sca \times n_\rho$). This is a significant saving compared with a brute-force computation. Starting from~\eqref{FIO}, $\alpha^\omega$ is written as a linear operator acting on the data $\Lambda^\omega$ through an integral kernel. {Following this more direct approach, the input, $\sf\Lambda^\omega$, will be flattened into a vector of size $n_\sca^2$, and the output, $\alpha^\omega$, will flattened into a vector of size $n_\eta n_\rho$, therefore the integral kernel is represented by a significantly larger matrix of size $(n_\eta n_\rho)\times n_\sca^2$.} This procedure not only ignores the fact that the kernel is independent of $\theta_j$, but also requires $n_\eta n_\sca$ times more unknowns, increasing the overall computational  cost.

\begin{figure}[h]
\centering
\includegraphics[width=\linewidth]{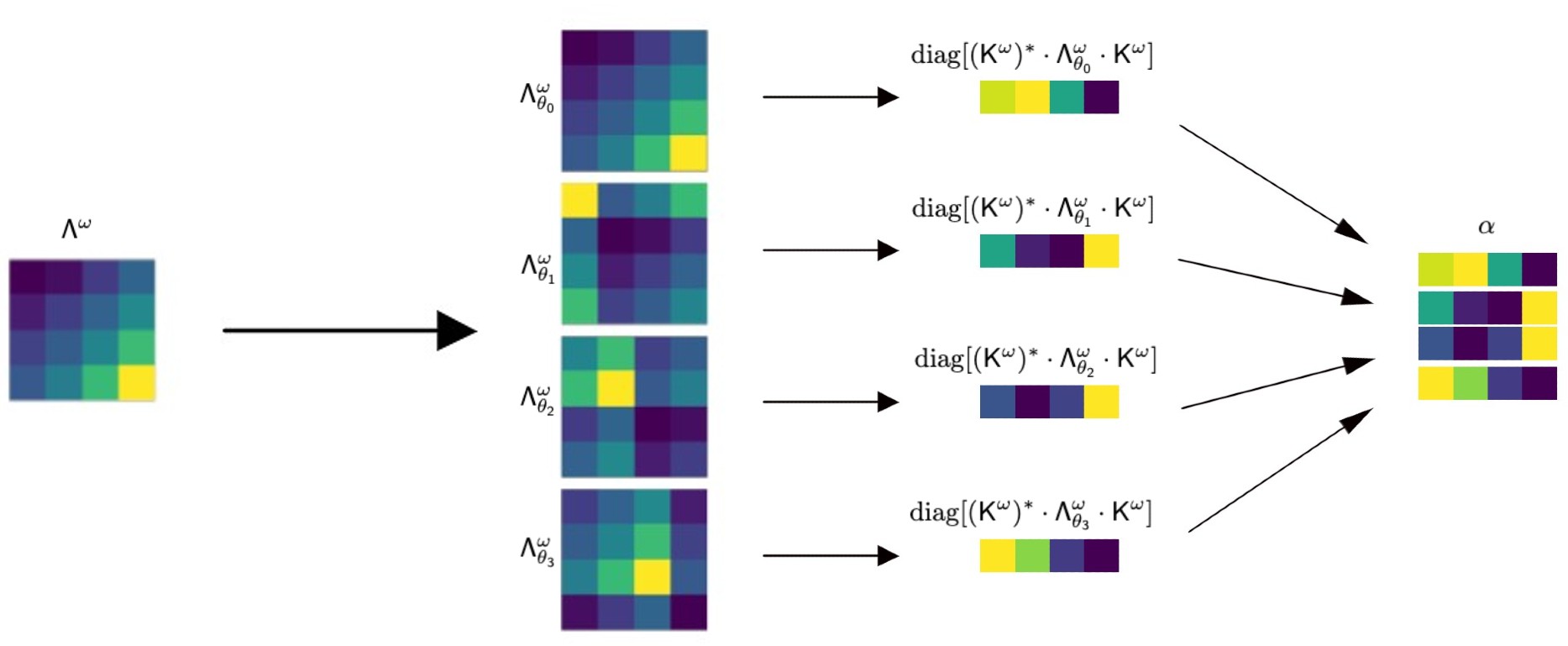}
\caption{The visualization of the application of the underlying equivariance. In the first step, the data matrix $\sf\Lambda^\omega$ is shifted to generate the other four $\sf\Lambda_{\theta_j}$ for $j=0,1,2,3$. Then, the Implementation II
$\sfx\mapsto\operatorname{diag}
[(\mathsf{K}^\omega)^\ast \cdot \sfx\cdot \mathsf{K}^\omega]$ is applied to all four $\sf\Lambda_{\theta_j}$, each of which outputs a row vector. Finally, they are concatenated to form the intermediate representation $\alpha^\omega$.}
\label{fig:visualization}
\end{figure}

In Figure~\ref{fig:visualization}, we visualize the application of the underlying equivariance on the back-scattering operator, using $n_{\sca}=n_\rho=n_\theta=4$ as an example for Implementation II.

In order to implement numerically the algorithm in Implementation I of~\eqref{eqn:two_forms}, we note that the kernel can be decomposed as following
\begin{equation}\label{eqn:KCS}
K^\omega =C^\omega-iS^\omega\,,\quad\text{with}\quad
C^\omega(\rho,t) = \cos(\rho\omega\cos(t))\quad\text{and}\quad S^\omega(\rho,t) = \sin(\rho\omega\cos(t))\,.
\end{equation}
In this case, we have
\begin{align*}
((F^{\omega})^* \Lambda^\omega)(\theta,\rho)
&= \iint_{[0,2\pi]^2} C^\omega(\rho,r)C^\omega(\rho,s)\Lambda^\omega_R(r+\theta,s+\theta)\,ds\,dr\,,\\
&+ \iint_{[0,2\pi]^2} S^\omega(\rho,r)S^\omega(\rho,s)\Lambda^\omega_R(r+\theta,s+\theta)\,ds\,dr\,,\\
&+ \iint_{[0,2\pi]^2} C^\omega(\rho,r)S^\omega(\rho,s)\Lambda^\omega_I(r+\theta,s+\theta)\,ds\,dr\,,\\
&- \iint_{[0,2\pi]^2} S^\omega(\rho,r)C^\omega(\rho,s)\Lambda^\omega_I(r+\theta,s+\theta)\,ds\,dr\,,
\end{align*}
where $\Lambda^{\omega}$ has been periodically extended outside the domain $(r,s)\in[0,2\pi]^2$ torus.

As such, the training of $\sfK^\omega$ can be translated to the training of $\sfC^\omega$ and $\sfS^\omega$, which are the discrete forms of $C^\omega$ and $S^\omega$ using the same discretization as $K^\omega$. Following the Implementation I in equation~\eqref{eqn:two_forms}, we find:
\begin{equation}\label{eqn:CS}
\alpha^\omega(\theta_j,\cdot) = \operatorname{ones}(1,n_\sca)\cdot
[\sfC^\omega \odot (  \sf\Lambda^\omega_{\theta_j,R}\cdot\sfC^\omega)+\sfS^\omega \odot (\sf\Lambda^\omega_{\theta_j,R}\cdot\sfS^\omega )+\sfC^\omega \odot (  \sf\Lambda^\omega_{\theta_j,I}\cdot\sfS^\omega)-\sfS^\omega \odot ( \sf\Lambda^\omega_{\theta_j,I} \cdot\sfC^\omega)]\,,
\end{equation}
where $\sf\Lambda^\omega_{\theta_j,R/I}$ are shifted real/imaginary parts of the data as defined in~\eqref{eqn:data_shift}. Clearly, the execution of $(\mathsf{F}^\omega)^\ast$ is written as the summation of four terms that share the {same sequence of operations. Namely,} each term is composed {of} one vector-matrix product, one Hadamard product, and one matrix-matrix product with the shifted data.

We summarize the implementation of~\eqref{eqn:CS}, in which the parameters $\sfC^\omega$ and $\sfS^\omega$ , instead of being constant matrices given by~\eqref{eqn:KCS}, are to be trained using data. Also, we replace the constant matrix $\operatorname{ones}(1,n_\sca)$ by trainable weights $\sfO^j$ for $j=1,2,3,4$ of the same size, so that we can get rid of the negative sign. We name the model as the uncompressed model.


\begin{algorithm}[H]
\begin{algorithmic}[1]
\caption{The application of $(F^\omega)^*$ in the uncompressed model}
\Statex Input: {$\sf\Lambda^\omega \in \mathbb{C}^{n_{sc} \times n_{sc}} $}
\Statex Output: {$(\sfF^\omega)^*$}{$\sf\Lambda^\omega \in \mathbb{C}^{n_{\eta} \times n_{\rho}}$}
\State Split the data in real and imaginary parts: $
\sf\Lambda^\omega=\Lambda^\omega_R+i\Lambda^\omega_I$
\For{$j<n_{sc}$}
     \State {$\alpha^\omega[j,:] \gets 
     \sfO^1 \cdot( \sfC\odot(\sf\Lambda^\omega_{\theta_j,R}\cdot \sfC))
   + \sfO^2 \cdot( \sfS \odot(\sf\Lambda^\omega_{\theta_j,R}\cdot \sfS))
   + \sfO^3 \cdot( \sfC\odot(\sf\Lambda^\omega_{\theta_j,I}\cdot \sfS)) 
   + \sfO^4 \cdot( \sfS\odot(\sf\Lambda^\omega_{\theta_j,I}\cdot \sfC))$}
\EndFor
\State \textbf{return} $\alpha^\omega$
\end{algorithmic}
\end{algorithm}
We note that the model presented in this section only utilizes the underlying equivariance. In the application of $(F^\omega)^\ast$, we represent the maps $S^\omega$ and $C^\omega$ directly as $n_\sca\times n_\sca$ matrices consisting of trainable weights and the size of $\bar{\Omega}$ scales logarithmically with respect to $n_\sca$ following a dyadic partitionning, which ensure that we capture high-frequency information while still taking advantage of the regularization power that low frequency data provides to the algorithm.  Hence, the number of trainable weights scales as $\mathcal{O}(n_\sca^2\log n_\sca)$, and the inference complexity is dominated by the computations of $n_\sca\times n_\sca$ matrices productions for a total of $4n_\sca\log n_\sca$ times, which scales $\mathcal{O}(n_\sca^4\log n_\sca)$.

Since the filtering operator is approximated by a 2-dimensional convolutional NN with constant-sized kernel, for which the size of filters and the number of layers scale at most linearly with respect to $n_\sca$, the number of trainable parameters scales $\mathcal{O}(n_\sca(\log n_\sca)^2)$, and the inference complexity is {$\mathcal{O}(n_\sca^3(\log n_\sca)^2)$}. It should be noted that the $\log n_\sca$ term comes from the dependence of the channel dimension on the number of frequencies. 

Therefore, for the uncompressed model, the total number of parameters scales $\mathcal{O}(n_\sca^2\log n_\sca)$ and the inference complexity is $\mathcal{O}(n_\sca^4\log n_\sca)$.

We should note that the equivariance in our setting is hardcoded into the neural-network design. The feature of equivariance is widely observed in many other scientific domains, such as computer vision, reinforcement learning, dynamics learning, or protein folding. This triggers the studies on preserving gauge invariance structure in its very general form, see for example~\cite{cohenc16:EquivariantCNN} and~\cite{bronstein2021geometric}. In the context of inverse scattering, we have found~\cite{FanYing:scattering} that also places focuses on preserving rotational equivariance. The compression, however, was not incorporated in their study.

\subsection{Application of the butterfly factorization}\label{sec:butterfly}
For the application of butterfly factorization, we adopt Implementation II in~\eqref{eqn:two_forms} of the discretized adjoint operator $(F^\omega)^\ast$, which results on the formula
\[\alpha^\omega(\theta_j,\cdot)=\sf((F^\omega)^*\Lambda^\omega)(\theta_j,\cdot) = \operatorname{diag}
[(\mathsf{K}^\omega)^\ast \cdot \sf\Lambda^\omega_{\theta_j}\cdot \mathsf{K}^\omega]\,.\]
Notice that in Section \ref{sec:propertybutterfly} we have shown that $\mathsf{K}^\omega$ admits a butterfly factorization, namely, it can be factorized as in~\eqref{eqn:butterfly}.
As a direct consequence, $(\mathsf{K}^\omega)^\ast \cdot \sf\Lambda^\omega_{\theta_j}\cdot \mathsf{K}^\omega$ is expanded as:
\begin{equation}\label{eqn:butterfly_2}
\sfV^L\sfH^{L-1}\cdots \sfH^{L/2}(\sfM^{L/2})^\ast(\sfG^{L/2})^*\cdots(\sfG^{L-1})^*(\sfU^L)^*\cdot {\sf\Lambda^\omega}\cdot \sfU^L\sfG^{L-1}\cdots \sfG^{L/2}\sfM^{L/2}(\sfH^{L/2})^*\cdots(\sfH^{L-1})^*(\sfV^L)^*\,.
\end{equation}

The form of~\eqref{eqn:butterfly_2} suggests the overall structure of the operator is composed of taking actions on the data $\sf\Lambda^\omega_R$ using $L+3$ layers:

\begin{itemize}
    \item Layer $U_{layer}(\sfx)$: In this layer, we apply the most internal action on the data: we sandwich the original data $\sf\Lambda^\omega_R$ by $\sfU^L$, namely:
    \[
    \sfx\to(\sfU^L)^*\cdot \sfx\cdot \sfU^L\,.
    \]
    \item Layer $G_{layer}(\sfx,\ell)$: there are $L/2$ $G_{layer}(\sfx,\ell)$ layers corresponding to $\ell =L/2\,,\cdots\,, L-1$, each of which sandwich the given data by the associated $\sfG^i$:
    \[
    \sfx \to (\sfG^\ell)^\ast\cdot\sfx\cdot\sfG^\ell\,,\quad \ell =L/2\,,\cdots\,, L-1\,.
    \]
    \item Layer SwitchResnet($\sfx$): In the middle layer, we are supposed to sandwich the given data by $\sfM^{L/2}$:
    \[
    \sfx \to (\sfM^{L/2})^\ast\cdot\sfx\cdot\sfM^{L/2}\,.
    \]
    Additionally, inspired by~\cite{MLZ}, we also lift the middle layer by integrating a Resnet structure into it. We call it a SwitchResNet. The rational will be explained in the description of the implementation later in this section.
    \item Layer $H_{layer}(\sfx,\ell)$: similar to the $G_{layer}(\sfx,\ell)$ layer. 
    \item Layer $V_{layer}(\sfx)$: similar to the $U_{layer}(\sfx)$ layer. 
\end{itemize}

We note that the application of $\sfK^\omega$ assumes small perturbation around $\eta_0$ in the linearized setting, and to lift it up to deal with nonlinear inverse scattering, we still would like to honor these symmetries. In computation, we keep the linearity of every action except lifting the application of $\sfM^{L/2}$ and its adjoint $(\sfM^{L/2})^\ast$ to the nonlinear setting and representing them by a ResNet. This is summarized in the pseudo-code in Algorithm~\ref{compressed}, see also~\cite{MLZ}. We name the corresponding model as the compressed model.

\begin{algorithm}[H]
\begin{algorithmic}[1]
\caption{The application of $(F^\omega)^*$ in the compressed model.}\label{compressed}
\Statex Input: {$\sf\Lambda^\omega \in \mathbb{C}^{n_{sc} \times n_{sc}} $}
\Statex Output: {$(\sfF^\omega)^*$}{$\sf\Lambda^\omega \in \mathbb{C}^{n_{\eta} \times n_{\rho}}$}
\State $\sfx\gets U_{layer}(\sf\Lambda^\omega)$
\For{$\ell$ in $\operatorname{range}((L-1), L/2-1, -1)$}:
\State $\sfx \gets G_{layer}(\sfx,\ell)$
\EndFor
\State $\sfx \gets \operatorname{SwitchResnet}(\sfx)$
\For{$\ell$ in $\operatorname{range}(L/2, L+1, +1)$}:
\State $\sfx \gets H_{layer}(\sfx,\ell)$ 
\EndFor
\State $ \sfx\gets V_{layer}(\sfx)$
\State \textbf{return} $\sfx$
\end{algorithmic}
\end{algorithm}

The reasoning for choosing everything linear but lifting the middle layer to SwitchResnet stems from the fact different layers are in charge of exchanging information at different levels. In particular, seen from Figure~\ref{fig:bf}, $\sfU^L$ and $\sfV^L$ are block-diagonal matrices and thus are automatically in charge of local information change. $\sfG^l$ and $\sfH^l$ are in charge of exchanging information across neighboring blocks and thus the action is also local. $\sfM^{L/2}$ finally is responsible for capturing the inherent non-locality of wave scattering. Since this action is taken on the condensed representation of the measured data, we choose to represent it using a non-linearly process.

There are multiple trade-offs in the implementation of different layers, and we glean through them below. As an example, we assume $\sf\Lambda^\omega_R$ is an $80\times 80$ matrix. As implied by the Implementation II of~\eqref{eqn:two_forms}, $(\sfK^\omega)^\ast$ (and similarly $\sfK^\omega$) is applied on each individual column (row) of $\sf\Lambda^\omega_R$ viewed as an independent vector, so below we describe the application of the network on one of the columns of the data matrix, which is an $80$ dimensional vector. To deal with this vector, we divide it into $2^4$ chunks with each chunk containing $5$ entries. As a result, calling $80=2^4\times 5$, we set $L=4$ and $s=5$.
\begin{itemize}
    \item[--]{\textbf{Implementing the $U_{layer}$ and $V_{layer}$ layers:}}\\
    As shown in Figure~\ref{fig:bf}, in the $U_{layer}$ and $V_{layer}$ layers, we sandwich the data by the block matrices $\sfU^L$ and $\sfV^L$. In the $U_{layer}$  layer, each block of $\sfU^L$, denoted as $\sfU^L_i$ with $i=1\,,\cdots 2^L=16$ extracts a local representation of the $i$-th section of the vector (an $s$-dimensional vector) and represent it by a $r$-dimensional vector. Used in our setting, each of $\sfU^L_{i}$, when applied to the $s=5$-dimensional vector (the $i$-th section of the vector), produces a $r$-dimensional vector. In total, the procedure generates a $2^Lr=16r$-dimensional vector. The number of trainable weights in the matrix is therefore $2sr2^L$ after separating the real and complex channels. Similarly for the $V_{layer}$ layer, each block of $\sfV$ transforms an $r$-dimensional local representation to a $s$-dimensional sampling points, which also requires $2sr2^L$ trainable weights.



    \item[--]{\textbf{Implementing the $G_{layer}$ and $H_{layer}$ layers:}}\\
    In the $G_{layer}(\sfx,\ell)$ layer, $(\sfG^\ell)^\ast$ and $\sfG^\ell$ are applied to $\sfx$, and in the $H_{layer}(\sfx,\ell)$ layer, $(\sfH^\ell)^\ast$ and $\sfH^\ell$ are applied to $\sfx$.
    Observing the Figure~\ref{fig:bf}, we see that in the $G_{layer}$ layer, two neighboring local representations are assimilated. Similarly, in the $H_{layer}$ layer, information in the representations will be decompressed by being splitted in two and locally redistributed.  To implement the $G_{layer}$ and $H_{layer}$ layers, we decompose the corresponding matrices into row or column transformations and block matrices, so that after the row or column transformations, the application of the block matrices can be implemented the same way as in $U_{layer}$ and $V_{layer}$ layers. The number of trainable weights in both of the layers is therefore $4r^22^L$.

    \item[--]{\textbf{Implementing the $M_{layer}$ layer:}}\\
    In the $M_{layer}$ layer, information are redistributed globally by the matrices $(\sfM^{L/2})^\ast$ and $\sfM^{L/2}$. The application of $M_{layer}$ layer is essentially the same as the $G_{layer}$ and $H_{layer}$ layers, where we apply the corresponding row or column transformations following by block matrices. Inspired by~\cite{MLZ}, we also integrate a non-linear module into the model by adding a Resnet of depth $n_{\mathsf{SR}}$ in the $M$ layer. Hence, we have a total of $2n_{\mathsf{SR}} r^22^L$ trainable weights.

\end{itemize}
Compared to the uncompressed model, the compressed model exhibits lower asymptotic scaling of the number of trainable parameters and a 
 lower inference complexity, albeit still super-linear. As explained at the end of Section~\ref{sec:equivariance}, the number of trainable parameters of the filtering operator scales $\mathcal{O}(n_\sca(\log n_\sca)^2)$ and the inference complexity of the filtering operator is {$\mathcal{O}(n_\sca^3(\log n_\sca)^2)$}. Nevertheless, for the application of $(F^\omega)^\ast$ in the compressed model, the number of parameters scales $(4sr2^L+8Lr^22^L+2n_{\mathsf{SR}}r^22^L)\log n_\sca=\mathcal{O}(r^2n_\sca(\log n_\sca)^2)$, and the inference complexity is $\mathcal{O}(r^2n_\sca^3(\log n_\sca)^2)$. 

Therefore, for the compressed model, the complexity of the number of trainable parameters is $\mathcal{O}(r^2n_\sca(\log n_\sca)^2)$ and the inference complexity is $\mathcal{O}(r^2n_\sca^3(\log n_\sca)^2)$.
In Table~\ref{tab:complexity}, we present the complexity of the trainable parameters and the inference time complexity for the compressed model and the uncompressed model. 

\begin{table}[h]
\begin{tabular}{ |p{3.5cm}||p{5cm}|p{5cm}| }
 \hline
 Complexity & Compressed model & Uncompressed model\\
 \hline
Parameters  & $ \mathcal{O}(r^2n_{\sca}(\log{n_{\sca}})^2)$
  &  $\mathcal{O}(n_{\sca}^2\log{n_{\sca}})$  \\
\hline
Inference time  &  $\mathcal{O}(r^2n_{\sca}^3 (\log{n_{\sca}})^2)$ &  $\mathcal{O}(n_{\sca}^4\log{n_{\sca}})$ \\
 \hline
\end{tabular}

\caption{Complexity of the number of trainable parameters and inference complexity between the uncompressed and the compressed models. The complexity of the uncompressed and compressed models was discussed in the Sections~\ref{sec:equivariance} and~\ref{sec:butterfly}, respectively.} 
\label{tab:complexity}
\end{table}

\section{Numerical Results}\label{sec:numerics}
In this section, we provide numerical evidences to showcase the capabilities of the current algorithm. In summary, what to be presented below suggest two findings:
\begin{itemize}
    \item Uncompressed model has the equivariance built in the formulation, leading to fewer trainable parameters. The amount of data needed for the training is consequently reduced. The performance of the numerical results are competitive compared to classical FWI and other NN models (such as wide-band butterfly network and Fourier neural operator) even with smaller amount of training points.
    \item Compressed model further reduces the number of trainable parameters by building in the butterfly structure. The performance of the model is biased towards media that are relatively smooth: The reconstruction of smooth media is accurate while the model losses the capability of capturing sub-Nyquist features. 
\end{itemize}

The training optimization formulation is presented in subsection~\ref{sec:optimization}. We then discuss data structure in subsection~\ref{sec:datasets}. These will be followed by subsection~\ref{sec:uncompressed} that discusses the uncompressed models. In subsections~\ref{sec:compressed}, we systematically compare reconstruction using various of means: The uncompressed model, the compressed model, the classical FWI and two other well-accepted NN models.

The data was generated using Matlab on a ARM-based MacBook Air (M1, 2020). The script usually took about 12 hours to generate a multi-frequency datasets of dimension $n_{\sca}=80$, of size $10000$, and with source frequencies 2.5, 5.0, and 10.0 Hz. The models presented in this paper were implemented using Tensorflow (2.4.1)~\cite{tensorflow2015} and ran on a PNY NVIDIA Quadro RTX 6000 graphics card.

\subsection{Optimization}\label{sec:optimization}

{As before, we denote our neural network as $\Phi_{\Theta}$ with parameters $\Theta$ and the set of frequencies of the data begin fed to the network by $\bar{\Omega}$, namely our network takes the form: 
\begin{equation} \eta = \Phi_{\Theta} (\{\Lambda^{\omega}\}_{\omega \in \bar{\Omega}})\,.
\end{equation}}

{For training we use tuples of media and scattering data, $(\eta^{[s]}, \{\Lambda^{\omega, [s]}\})$, where $[s]$ is the index for the training set. We point out that different problems will have different $\bar{\Omega}$ depending of the target resolution, also the network will change depending on whether the operators are compressed.}

{We train the networks by minimizing the mean square error between the network produced media and the groundtruth media (the media used to generate the input data), i.e.,
\begin{equation}
\min_{\Theta} \frac{1}{N_s}\sum_{s =1}^{N_s} \| \Phi_{\Theta} (\{\Lambda^{\omega, [s]}\}_{\omega \in \bar{\Omega}}) - \eta^{[s]} \|^2\,.
\end{equation}  
We use Adam optimizer with learning rate chosen as $3\times10^{-4}$, and batch size as $16$ together with an exponential scheduler. The learning scheduler was set as Tensorflow's~\cite{tensorflow2015} \texttt{ExponentialDecay} with a decay rate of 0.96 after every 50 plateaus steps with \texttt{staircase} option set as true. We chose the Adam optimizer~\cite{kingma2015adam} and we train the model after 100 epochs.  The trainable weights were all randomly initialized with the glorot\_uniform distribution~\cite{glorot}.}
{For comparison we report the relative root mean square error (RMSE) given by 
\begin{equation}
\frac{1}{N_t}\sum_{s =1}^{N_t} \frac{\| \Phi_{\Theta} (\{\Lambda^{\omega, [s]}\}_{\omega \in \bar{\Omega}}) - \eta^{[s]} \|}{ \| \eta^{[s]}\|} \,,
\end{equation}
where $N_t$ is the size of the testing set. }


\subsection{Datasets}\label{sec:datasets}

The datasets consist of media and corresponding wide-band far-field patterns {at} three different frequencies, which are rescaled appropriately depending on the target resolution. The far-field patterns as the data were generated by solving for equation~\eqref{eqn:scattereqn} using numerical finite difference method in Matlab. 
The computational domain was $[-0.5,0.5]^2$ discretized with a equispaced mesh of $80 \times 80$ points for medias of resolution $80 \times 80$ pixels ($n_\eta=80$). The radiation boundary conditions was implemented using the perfectly matched layer (PML)\cite{Berenger:PML} with order 2 and intensity 80. The wide-band data was sampled with a homogeneous background wave field at source frequencies 2.5, 5, and 10 Hz, see Section~\ref{sec:preliminary}, for which the effective wavelength is 8 points per wavelength (PPW). In particular, we use $n_\sca=80$ receivers and sources, where receivers geometry are aligned with the directions of sources, i.e. 80 equiangular directions. For media of different resolutions, we generate their far-field patterns by sampling at proportionally scaled frequencies.  

\begin{figure}[h]
\center
\includegraphics[width=1\textwidth]{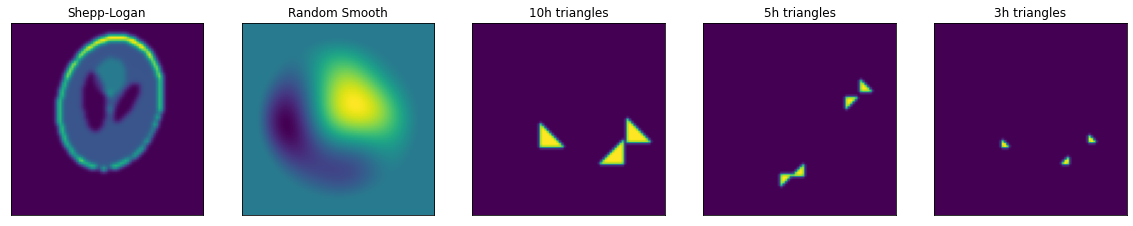}
\caption{Examples of the five different media used to bechmarkthe model (from left to right) Shepp-Logan phantom, smooth random perturbation, and the rotating triangles of different sizes.}\label{fig:media}
\end{figure}

In our experiments we use 5 different categories of media: 
\begin{itemize}
\item The well-known Shepp-Logan phantom,  which was created in 1974 by Larry Shepp and Benjamin F. Logan to represent a human head~\cite{SL}. The medium has a strong discontinuity modeling an uneven skull, which produced a strong reflection, which in return renders the recovery of the interior features challenging for classical methods. 
\item Random smooth perturbations, which are generated by smoothing out some randomly distributed points of some random values by a Gaussian kernel.  They are used to study the behavior of the algorithm in the case of diffraction.
\item Media consisting of triangles of different sizes, in particular, triangles of size $3$, $5$ and $10$ number of pixels, which are randomly located and oriented, and it is possible for them to overlap with each other. In this case we test the capacity of the algorithm to image consisting of small scatterers that are slightly below  sub-Nyquist in size. For clarity we name the dataset following the sizes of the triangles, e.g., `10h triangles' are composed of triangles with side length of 10 pixels.  
\end{itemize}
In Figure~\ref{fig:media}, we provide one example for each of the five category. The training and the test are all conducted within one category.

\subsection{Uncompressed Model}\label{sec:uncompressed}

The uncompressed model has the equivariance built in. As a consequence, only a small number of parameters need to be trained, which in return requires fewer numbers of samples for the training. In what follows we demonstrate a few properties of uncompressed model. In particular,
\begin{enumerate}
    \item \textbf{Low sample complexity}: the validation error quickly saturates as we increase the number of training points;
    \item \textbf{Improved stability with wide-band data}: similarly to the classical FWI method, which produces better reconstruction with wide-band data, the uncompressed model provides better reconstruction when data at multiple frequencies are provided;
    \item \textbf{Improved accuracy with higher frequency data}: when the models are trained on media of finer resolution and scattered data of higher dimension that are generated by probing waves of higher frequencies, it produces results with lower relative validation error.  In particular, it recovers the fine-grained details better.
\end{enumerate}

\begin{figure}[h]

\includegraphics[width=1.0\textwidth]{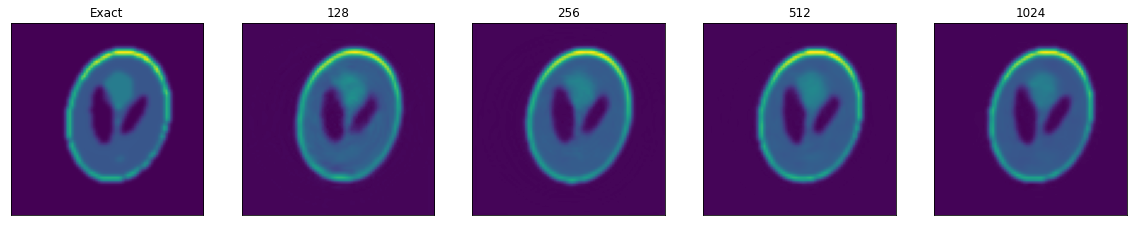}
\includegraphics[width=1.0\textwidth]{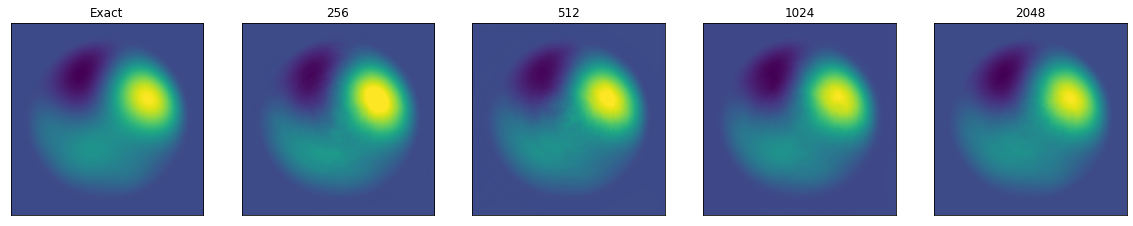}
\includegraphics[trim={0 0 0 1cm},clip,width=1.0\textwidth]{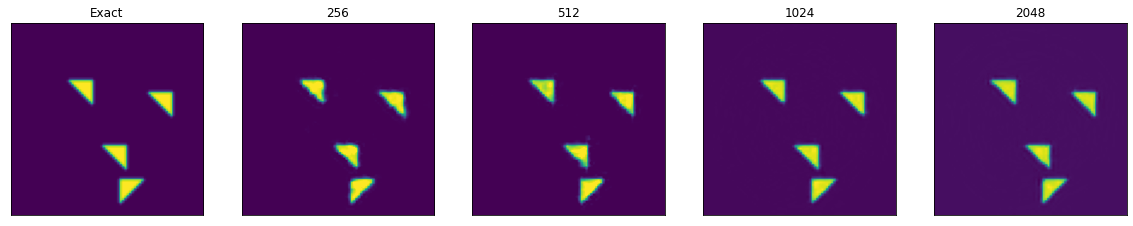}
\includegraphics[trim={0 0 0 1cm},clip,width=1.0\textwidth]{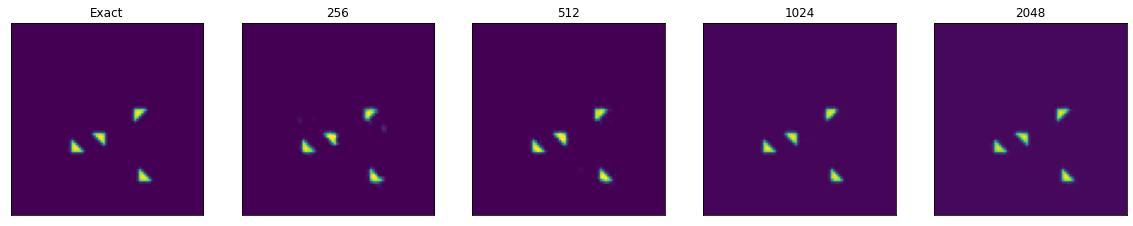}
\includegraphics[trim={0 0 0 1cm},clip,width=1.0\textwidth]{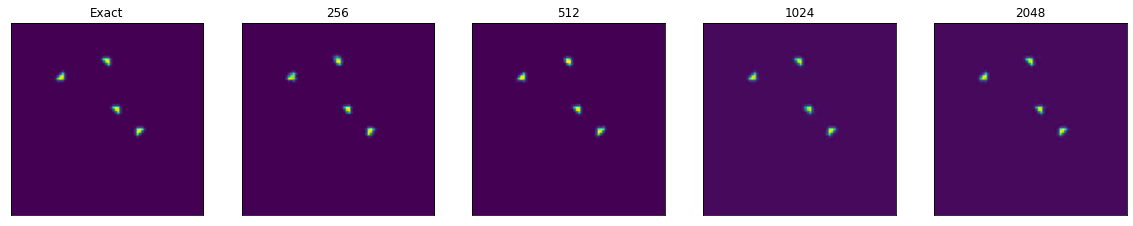}

\caption{Illustration of the numerical performance of the Uncompressed model. Each row presents results for a different media using different numbers of training points in the training process. The exact media is shown in the left-most column. The reconstruction utilizes the wide-band frequency datasets.}\label{fig:illustration_uncompressed}

\end{figure}

In Table~\ref{samplecomplexity}, we show how the relative validation error depends on the number of training data points, and probing wave frequency for different media. The scattered field is generated using $n_\sca = 80$ and the media in a domain is discretized with an equispaced mesh of $n_\eta = 80$. The number of trainable parameters for reconstruction with one frequency is 46530 and that for reconstruction with three frequencies is 88186. {Table~\ref{samplecomplexity} shows that} using wide-band data consistently produces more accurate reconstructions; the equivariance property helps to drastically reduce the number of training points (the accuracy has already stagnated using a few thousands training points).
{Figure~\ref{fig:illustration_uncompressed} illustrates one representative instance of the reconstruction for each of the different media mentioned above, using different numbers of training points. In this case, the model uses wide-band frequency datasets.}

\begin{table}[h]
\begin{tabular}{ |p{3cm}||p{3cm}|p{3cm}|p{3cm}|p{3cm}|  }
 \hline
 \multicolumn{5}{|c|}{Shepp–Logan Phantom} \\
 \hline
 \#Sample $\backslash$ Frequency & 2.5 Hz & 5 Hz & 10 Hz & 2.5 \& 5 \& 10 Hz\\
 \hline
 64 & 12.430 \% & 10.119 \% & 15.938 \% & 13.236\%  \\
128 & 10.619 \% & 7.919 \% & 11.403 \%  & 8.293 \%  \\
256 & 10.334 \% & 6.754 \% & 9.798 \%   & 6.503 \% \\
512 & 10.077 \% & 6.935 \% & 7.954 \%   & 5.124 \% \\
1024 & 10.774 \% & 6.721 \% & 8.125 \%  & 5.306 \% \\
 \hline
\end{tabular}

\begin{tabular}{ |p{3cm}||p{3cm}|p{3cm}|p{3cm}|p{3cm}|  }
 \hline
 \multicolumn{5}{|c|}{Random Smooth Perturbation} \\
 \hline
 \#Sample $\backslash$ Frequency & 2.5 Hz & 5 Hz & 10 Hz & 2.5 \& 5 \& 10 Hz\\
 \hline
128 & 6.446 \% & 10.323 \% & 22.811 \%  &  7.182 \%  \\
256 & 5.240 \% & 7.649 \% & 21.127 \%   &  5.456 \%  \\
512 & 4.490 \% & 7.386 \% & 11.681 \%   &  4.898 \%  \\
1024 & 4.385 \% & 9.029 \% & 10.641 \%  &  3.957 \%  \\
2048 & 5.386 \% & 8.237  \% & 11.403  \%  &  4.102 \%  \\

 \hline
\end{tabular}

\begin{tabular}{ |p{3cm}||p{3cm}|p{3cm}|p{3cm}|p{3cm}|  }
 \hline
 \multicolumn{5}{|c|}{10h triangles} \\
 \hline
 \#Sample $\backslash$ Frequency & 2.5 Hz & 5 Hz & 10 Hz & 2.5 \& 5 \& 10 Hz\\
 \hline

256  & 47.911  \% & 18.432   \% & 14.858  \%  &  15.701  \%  \\
512  & 44.270  \% & 14.010  \% &  11.066 \%  &  9.012  \%  \\
1024 & 40.767  \% & 11.777 \% &   9.151 \%  &  7.777  \%  \\
2048 & 41.380  \% & 11.551  \% &  8.297  \%  &  7.068  \%  \\

 \hline
\end{tabular}

\begin{tabular}{ |p{3cm}||p{3cm}|p{3cm}|p{3cm}|p{3cm}|  }
 \hline
 \multicolumn{5}{|c|}{5h triangles} \\
 \hline
 \#Sample $\backslash$ Frequency & 2.5 Hz & 5 Hz & 10 Hz & 2.5 \& 5 \& 10 Hz\\
 \hline

256  &  59.616  \% & 37.625  \% &  11.020  \%  &   12.611   \%  \\
512  &  50.114  \% & 22.679  \% &  8.680  \%  &   8.842  \%  \\
1024 &  49.184  \% & 16.892  \% &  6.963  \%  &   6.061  \%  \\
2048 &  47.257  \% & 16.036  \% &  6.004  \%  &   6.230  \%  \\

 \hline
\end{tabular}

\begin{tabular}{ |p{3cm}||p{3cm}|p{3cm}|p{3cm}|p{3cm}|  }
 \hline
 \multicolumn{5}{|c|}{3h triangles} \\
 \hline
 \#Sample $\backslash$ Frequency & 2.5 Hz & 5 Hz & 10 Hz & 2.5 \& 5 \& 10 Hz\\
 \hline

256  & 59.327  \% & 33.156  \% &  12.054  \%  &   12.774  \%  \\
512  & 56.135  \% & 39.521  \% &  7.499    \%  &  10.209  \%  \\
1024 & 42.683  \% & 24.825  \% &  6.864   \%  &   6.597  \%  \\
2048 & 43.977  \% & 19.458  \% &  6.472   \%  &   5.902  \%  \\
 \hline
 \end{tabular}

 \caption{Relative {root mean square} test error of reconstruction given by the equivariant uncompressed model for difference media, using data at different frequencies, and with different sizes of training set. Each experiment consisted of \#Sample training points as recorded in the leftmost column and was tested against an independent testing dataset with 200 points.  Each table denotes a set of experiments with the type of media indicated by the table's name. Experiments on the first three column use monochromatic data sampled at source frequencies 2.5, 5, or 10 Hz, and experiments on the last column use wide-band data sampled consisting of data sampled at all three source frequencies. 
}\label{samplecomplexity}
\end{table}

Moreover, we also observe in Table~\ref{samplecomplexity} the unstable behavior of the reconstruction for the random smooth perturbations, particularly in the diffractive regime, in the sense that the error of reconstruction increases as the frequency increases for a single frequency, which is often a sign of cycle-skipping. This is solved by using wide-band data. The model is able to super-resolve some simple geometrical figures, such as the small rotating triangles, which is often a challenging problem using classical optimization/signal processing tools. We are able to reconstruct triangles of size 3 pixels, even if the wavelength is around 8 pixels.

\begin{figure}[h]
\center
\includegraphics[width=1\textwidth]{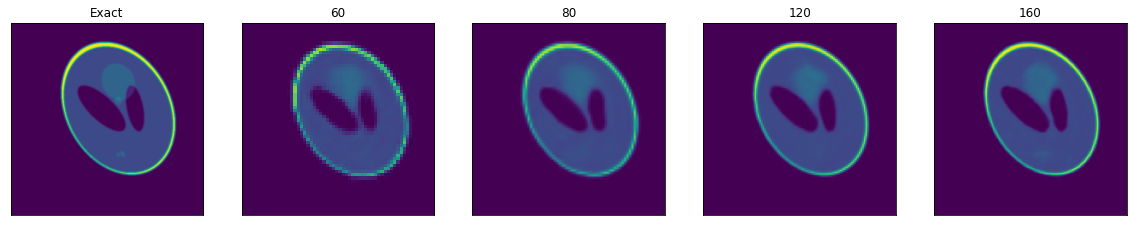}
\caption{Instances of the reconstruction of the Shepp-Logan phantom at different resolutions using the uncompressed method. In particular, the back-scattering operation of the method has more trainable weights catering to the larger input size, and the filtering operation, approximated by the 2D convolutional NN, has the same kernel size, number of filters and depth for all experiments. The left to right: the original Shepp-Logan phantom at its native resolution $n_\eta = 480$, the other four pictures are reconstructions of the downsampled media of size: $n_\eta=60,80,120$ and $160$ 
using training points of corresponding dimension: $n_\sca=60,80,120$ and $160$, repectively. The training points of higher resolutions are generated with probing wave of higher frequencies: $\frac{2.5}{80} n_\sca, \frac{5}{80} n_\sca$, and $\frac{10}{80} n_\sca$ Hz.}\label{fig:resolution}
\end{figure}

We also investigate how the number of trainable parameters scale according to the dimension of the data ($n_\sca$). Richer data is expected to produce better reconstruction. In our setting, as the data increases its dimension ($n_\sca$ becomes big), we expect the model to produce better results. To do so, we first randomly generate $1024$ samples of media at a native resolution with $n_\eta = 480$ and downsample them to a coarser resolution using $n_{\eta} = 60, 80, 120 \text{ and } 160$. At each coarse level of the media, we view the media as the reference and generate the wide-band far-field patterns using $n_{\sca}=n_{\eta}$, so that the frequencies of the probing waves scale proportionally to $n_{\sca}$. This builds the dataset for the training. Upon validation, one projects the reconstruction back to the native resolution at $n_\eta=480$. In Table~\ref{tab:size_and_rel_err}, we present at different coarse resolution, the number of trainable parameters, the average relative validation error at the training resolution level, and the relative validation error when the reconstruction gets interpolated back to the finest level (termed Rel err native res). In Figure~\ref{fig:resolution} we plot one example of the reconstruction at different resolution levels.

\begin{table}[H]
\begin{tabular}{ |p{3.5cm}||p{2cm}|p{3cm}|p{2.5cm}|p{3cm}|  }
 \hline
 \multicolumn{5}{|c|}{Shepp-Logan} \\
 \hline
 Resolution $\backslash$ Attributes & \#Parameters & Av inference time (s) & Rel err train res & Rel err native res\\
 \hline

60  & 70906  &  0.015  &  9.623  \%  &  27.719  \%  \\
80  & 88186  &  0.030  &  7.286  \%  &   22.854 \%  \\
120 & 137146  & 0.090 &  6.370   \%  &  11.192  \%  \\
160 & 205306  & 0.264 &   5.586  \%  &   7.223  \%  \\

 \hline
\end{tabular}
\caption{Statistics of the equivariant models for different resolutions of the Shepp-Logan dataset. The table shows the number of trainable parameters and the average inference time for the uncompressed model for Shepp-Logan phantom training points of resolution $n_\eta=60, 80, 120,$ and $160$. The last two columns records relative errors at the training resolution and at the native resolution for each test.}\label{tab:size_and_rel_err}
\end{table}

It is clear that if the datasets are prepared at a finer resolution, the number of trainable parameters increases accordingly, and so does the training time. The validation errors, both the relative error on the training resolution, and the relative error when projected back to the native resolution, also decrease.

\subsection{Validation of the Rotational Equivariance}

\begin{figure}[h]
\center
\includegraphics[width=1\textwidth]{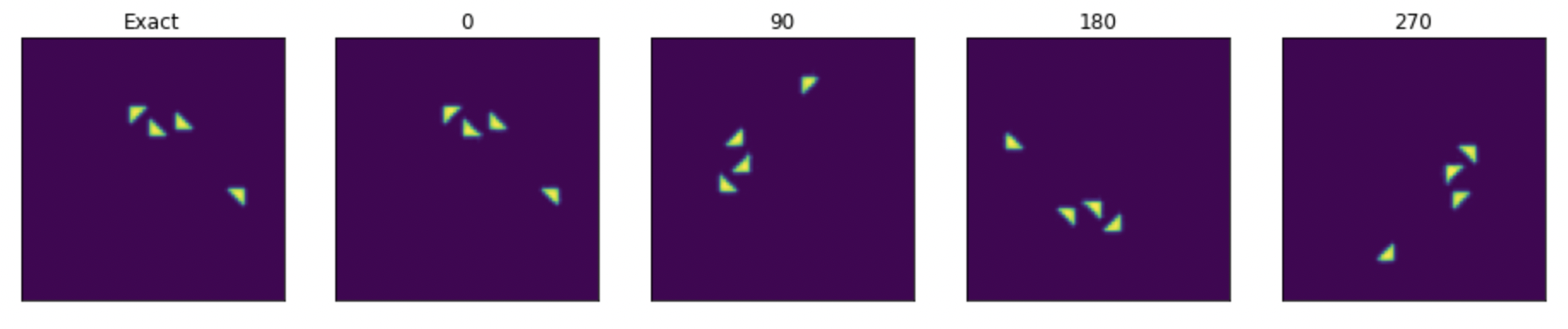}
\caption{Instances of the reconstruction of the 5h triangles rotated at 0, 90, 180, and 270 degrees.}\label{fig:rotation}
\end{figure}

We demonstrate the rotational equivariance of the back-scattering operator and the translational equivariance of the filtering operator, which are hard-wired into the neural network. In particular, we create four datasets by manually rotate the entire training set by a certain degree (0, 90, 180 and 270 degrees in four experiments), we then train four models: each corresponding to each new dataset. We consider as a base training set, the set of 5h triangles containing $2048$ data pairs. After training, we test each model on $500$ data points for validation. Table~\ref{tab:rotaionerror}, shows the average relative validation error for each rotated dataset. We can observe that regardless of the training data, the mean validation error of the reconstruction remains constant. In addition, Figure~\ref{fig:rotation}, shows one testing example of the reconstruction at different rotations: the reconstructed result is universally accurate across different rotations.

\begin{table}[h]
\center
\begin{tabular}{ |p{2.5cm}||p{2.5cm}|  }
 \hline
 Rotation (Degree) & Validation Error \\

 \hline 
$0 ^{\circ}$    & 5.599 \% \\
$90^{\circ}$   & 5.586  \%  \\
$180^{\circ}$   & 5.573  \%  \\
$270^{\circ}$   & 5.588 \%  \\
 \hline
\end{tabular}
\caption{The validation errors of the neural network model for datasets rotated at 0, 90, 180, and 270 degrees.}\label{tab:rotaionerror}
\end{table}

\subsection{Compressed Model}\label{sec:compressed}
In this section we present the numerical results produced by the compressed model, for which the numerical performance of the uncompressed model in Section \ref{sec:uncompressed} established a baseline. Moreover, we also include the numerical results from two widely accepted NN models: the wide-band butterfly network (WBBN) from~\cite{MLZ}, and the Fourier neural operator (FNO) from~\cite{fno} to form a comparison baseline with other ML methods. {A brief presentation of WBBN and FNO can be found in the appendices.} Since the compressed model calls for the construction of the butterfly factorization, the number of trainable parameters is further reduced. 

Throughout our simulation, we use training points of resolution at $n_\eta = 80$. This translates to, using the notation from Section~\ref{sec:butterfly}, $L=4$ and $s=5$. To train the uncompressed model, the compressed model, the WBBN and the FNO we use datasets consisting of far-field pattern data of dimension $n_\sca = 80$. The training set contains 2048 data points and the test set consists 500 data points. The batch size is chosen to be 16 for all the tests.

\begin{figure}[h]
\includegraphics[width=1.0\textwidth]{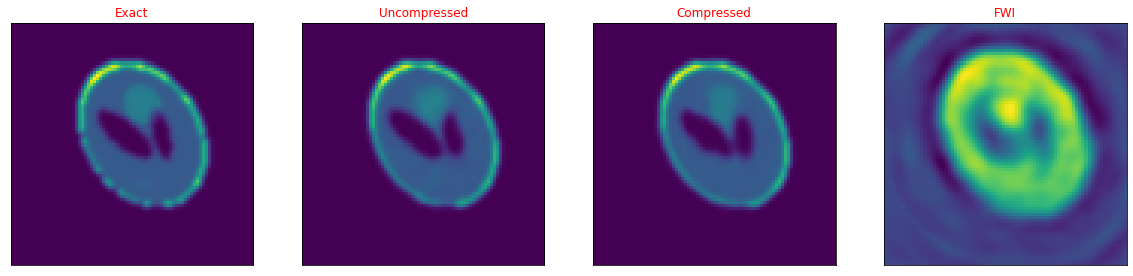}
\includegraphics[trim={0 0 0 1cm},clip,width=1.0\textwidth]{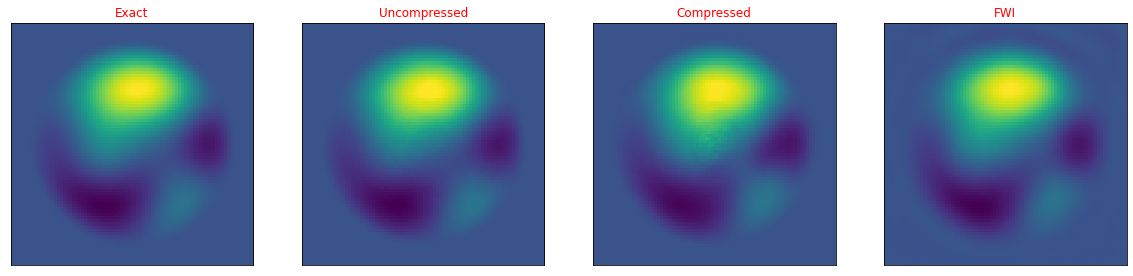}
\includegraphics[trim={0 0 0 1cm},clip,width=1.0\textwidth]{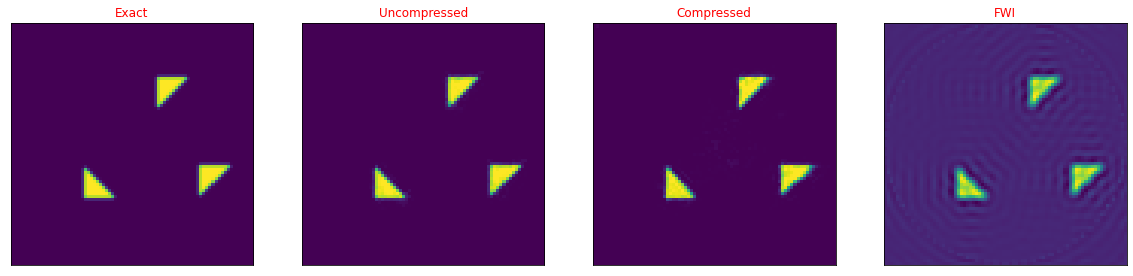}
\includegraphics[trim={0 0 0 1cm},clip,width=1.0\textwidth]{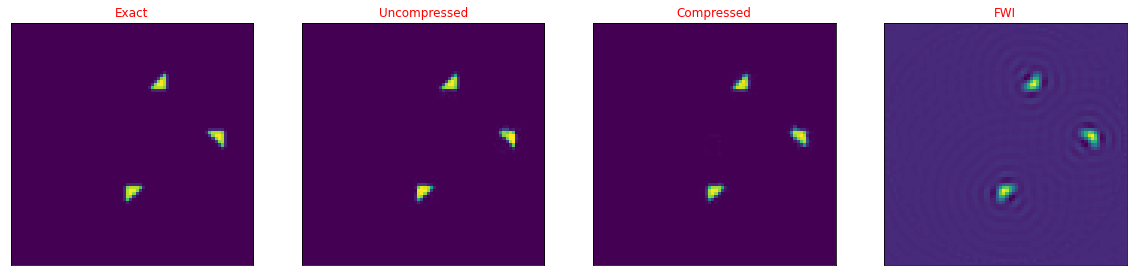}
\includegraphics[trim={0 0 0 1cm},clip,width=1.0\textwidth]{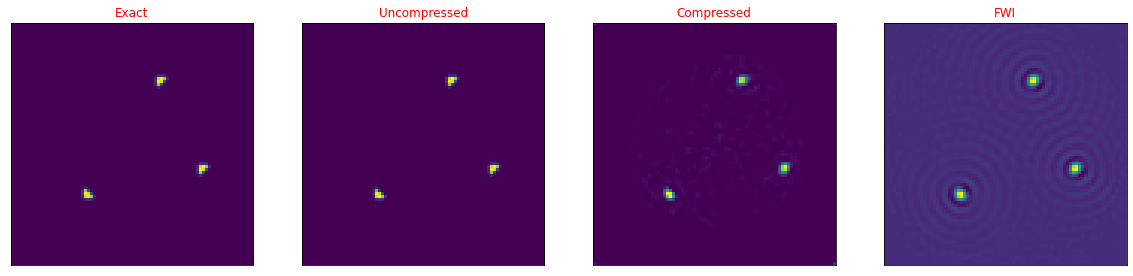}

\caption{Comparisons of reconstructions of various media by the uncompressed model, the compressed model, and the classical FWI. The media used are the Shepp-Logan phantom (the first row), the random smooth perturbations (the second row), and rotating triangles of different sizes. For the uncompressed and the compressed models, the dimension of data is chosen as $n_\sca=80$, and the number of data points is 2048.}\label{fig:comparison_FWI}
\end{figure}

In Figure~\ref{fig:comparison_FWI} we compare the uncompressed model result, compressed model result with the classical FWI that is performed by running the optimization with the lowest frequency first and using the previous results for the optimizations with the higher frequencies iteratively. It is clear that:
\begin{itemize}
    \item[1.] \textbf{Out-performance over FWI:} Both the uncompressed model and the compressed model achieve higher accuracy on all five kinds of media than the classical FWI;
    \item[2.] \textbf{Loss of super-resolution:} The compressed model has difficulty capturing super-resolution. It performs relatively well for 5h triangles, but lost the super-resolution for 3h triangle cases.
\end{itemize}

{For the WBBN model from~\cite{MLZ}, we test the model using the same hyperparameters as chosen in the paper. 
The initial learning rate was set as $5\times 10^{-3}$ and the scheduler was set as Tensorflow's~\cite{tensorflow2015} ExponentialDecay with a decay rate of 0.95 after every 2000 plateaus steps with staircase set to true. Adam optimizer~\cite{kingma2015adam} is employed and we terminate training after 150 epochs. Additionally, all trainable weights were randomly initialized with glorot\_uniform.} 

{To train FNO, the initial learning rate was set to be $10^{-3}$ and the scheduler was set as Pytorch's~\cite{pytorch} \texttt{StepLR} module with decay rate of 0.5 after every 100 plateaus steps with the \texttt{staircase} option set to true. Similarly, Adam optimizer~\cite{kingma2015adam} is employed and the training is terminated after 100 epochs.}

{The number of trainable parameters and the validation errors are listed in Table~\ref{tab:comparison_NNs}. In the case of Shepp-Logan phantom and the smooth perturbations, the performance of all four NN models is comparable even though the numbers of trainable parameters are drastically smaller for the compressed and uncompressed models, the two models proposed in the current paper.}

{However, for the small triangles that are below sub-Nyquist in size, both WBBN and FNO produce} {considerably worse results than our proposed method}. It should be noted that, for both WBBN and FNO, the numbers of trainable parameters are significantly greater than those of the uncompressed model and the compressed model,  which possibly results in over-fitting of the model when trained with a dataset of size 2048. Indeed, when trained with larger datasets, the WBBN and the FNO produce much better results. In the Table~\ref{tab:relerr_20000}, we list the validation errors of the WBBN and FNO with trained with 20000 data points using 10h, 5h, and 3h triangles.

\begin{table}[h]
\center
\begin{tabular}{ |p{2.5cm}||p{2.5cm}|p{2.5cm}|p{2.5cm}|p{2.5cm}|  }
 \hline
 Media $\backslash$ Model & Uncompressed & Compressed & WBBN & FNO \\
 \hline
\#parameters & 88,186 & 73,210 & 1,914,061 & 1,188,385 \\
 \hline 
Shepp-Logan  &   5.124 \%   &   8.306 
\%  & 9.819   \%  & 7.981 \%  \\
Random smooth &  3.957 \%   &   6.866  
\% & 8.268  \%  & 4.289 \%  \\
10h triangles &  7.068 \%    &   15.751 
\%  &  42.545  \%  &  42.633 \%  \\
5h triangles &   6.061 \%   &    16.144 
\%  & Untrainable   & 47.692 \%  \\
3h triangles &  5.902  \%   &   38.651 
\% & Untrainable   & 49.348 \%  \\
 \hline
\end{tabular}
\caption{Comparison of the {relative RMSE} for various NN models in the five different media categories.  The media reconstructed have resolution $n_\eta=80$. The data is generated using probing wave of frequencies $2.5$, $5$, and $10$ Hz, and the number of data points in the training set is $2048$. {Notice that the numbers of trainable parameters of the compressed model and the uncompressed model only have a small difference because of a relatively large constant in the complexity of the latter.} }\label{tab:comparison_NNs}
\end{table}

\begin{table}[h]
\center
\begin{tabular}{ |p{2.5cm}||p{2.5cm}|p{2.5cm}|  }
 \hline
 Media $\backslash$ Model & WBBN & FNO \\

 \hline 
Random smooth & 1.712 \%   & 2.253 \% \\
10h triangles & 2.858  \%  & 8.588  \%  \\
5h triangles &  3.521  \%  & 30.004  \%  \\
3h triangles &  5.119  \%  & 14.675 \%  \\
 \hline
\end{tabular}
\caption{The validation errors of the WBBN and FNO in the three media that contains small scattered are slightly below sub-Nyquist in size. The number of the data points is 20000.}\label{tab:relerr_20000}
\end{table}

\section{Conclusion}
In this manuscript, we perform numerical study on using various neuron networks to solve the 2D inverse scattering problem. In particular, we propose to build equivariance properties and the butterfly structures in the NN architecture design. By leveraging the underlying equivariance of the problem, we propose a model (termed the uncompressed model) that respects the equivariance of the back-scatter operator and reduces the number of trainable parameters, and by incorporating butterfly expansion structure we propose the second model (termed the compressed model) that further reduces the number of trainable parameters. The two models both outperform the classical FWI. For smooth media, they perform equally well as other NN architecures with fewer trainable parameters and smaller training sets, but the compressed model lost the super-resolution on sub-Nyquist features during the compression procedure.



\section{Acknowledgement}
The work of Q.\textasciitilde L. and B.\textasciitilde Z. is supported in part by the UW-Madison Data Initiative and the Office of Naval Research under the grant ONR-N00014-21-1-2140. The work of L.\textasciitilde Z.-N. is supported in part by the National Science Foundation under the grant DMS-2012292. In addition, Q.\textasciitilde L. and L.\textasciitilde Z.-N. are supported by the NSF TRIPODS award 1740707. The views expressed in the article do not necessarily represent the views of the any funding agencies. The authors are grateful for the support.

\bibliographystyle{unsrt} 

\bibliography{references} 

\begin{thebibliography}{10}

\bibitem{radar}
Brett Borden.
\newblock Mathematical problems in radar inverse scattering.
\newblock {\em Inverse Problems}, 18:R1 -- R28, 2001.

\bibitem{sonar}
C.~H. Greene, P.~H. Wiebe, J.~Burczynski, and M.~J. Youngbluth.
\newblock Acoustical detection of high-density krill demersal layers in the
  submarine canyons off georges bank.
\newblock {\em Science}, 241(4863):359--361, 1988.

\bibitem{seismic}
Arthur~B Weglein, Fernanda~V Ara{\'{u}}jo, Paulo~M Carvalho, Robert~H Stolt,
  Kenneth~H Matson, Richard~T Coates, Dennis Corrigan, Douglas~J Foster,
  Simon~A Shaw, and Haiyan Zhang.
\newblock Inverse scattering series and seismic exploration.
\newblock {\em Inverse Problems}, 19(6):R27--R83, oct 2003.

\bibitem{geophysics}
Dirk~J. Verschuur and A.~J. Berkhout.
\newblock Estimation of multiple scattering by iterative inversion; part ii,
  practical aspects and examples.
\newblock {\em Geophysics}, 62:1596--1611, 1997.

\bibitem{medical}
Tommy Henriksson, N.~Joachimowicz, Christophe Conessa, and Jean-Charles
  Bolomey.
\newblock Quantitative microwave imaging for breast cancer detection using a
  planar 2.45 ghz system.
\newblock {\em Instrumentation and Measurement, IEEE Transactions on}, 59:2691
  -- 2699, 11 2010.

\bibitem{Kirsch}
Andreas Kirsch.
\newblock {\em An Introduction to the Mathematical Theory of Inverse Problems}.
\newblock 2021.

\bibitem{garnier2016passive}
J.~Garnier and G.~Papanicolaou.
\newblock {\em Passive Imaging with Ambient Noise}.
\newblock Cambridge Monographs on Applied and Computational Mathematic.
  Cambridge University Press, 2016.

\bibitem{Hahner}
Peter H\"{a}hner and Thorsten Hohage.
\newblock New stability estimates for the inverse acoustic inhomogeneous medium
  problem and applications.
\newblock {\em SIAM Journal on Mathematical Analysis}, 33(3):670--685, 2001.

\bibitem{ChenDingLiZepeda}
Shi Chen, Zhiyan Ding, Qin Li, and Leonardo Zepeda-Núñez.
\newblock High-frequency limit of the inverse scattering problem: asymptotic
  convergence from inverse helmholtz to inverse liouville, 2022.

\bibitem{whittaker_1915}
E.~T. Whittaker.
\newblock Xviii.—on the functions which are represented by the expansions of
  the interpolation-theory.
\newblock {\em Proceedings of the Royal Society of Edinburgh}, 35:181–194,
  1915.

\bibitem{shanon1948:criterion}
C.~E. Shannon.
\newblock A mathematical theory of communication.
\newblock {\em The Bell System Technical Journal}, 27(3):379--423, 1948.

\bibitem{Courant1928berDP}
Richard Courant, K.~Friedrichs, and Hans Lewy.
\newblock {\"U}ber die partiellen differenzengleichungen der mathematischen
  physik.
\newblock {\em Mathematische Annalen}, 100:32--74, 1928.

\bibitem{Symes_Chen_Minkoff:2020}
William~W. Symes, Huiyi Chen, and Susan~E. Minkoff.
\newblock {\em Full-waveform inversion by source extension: Why it works},
  pages 765--769.
\newblock 2020.

\bibitem{ZepedaDemanet:the_method_of_polarized_traces}
L.~Zepeda-N\'u\~nez and L.~Demanet.
\newblock The method of polarized traces for the {2D} {H}elmholtz equation.
\newblock {\em J. Comput. Phys.}, 308:347 -- 388, 2016.

\bibitem{EngquistYing:Sweeping_PML}
B.~Engquist and L.~Ying.
\newblock Sweeping preconditioner for the {H}elmholtz equation: moving
  perfectly matched layers.
\newblock {\em Multiscale Model. Sim.}, 9(2):686--710, 2011.

\bibitem{linearrecursion}
Carlos Borges, Adrianna Gillman, and Leslie Greengard.
\newblock High resolution inverse scattering in two dimensions using recursive
  linearization, 2016.

\bibitem{fwi}
Yunyue~Elita Li and Laurent Demanet.
\newblock Full waveform inversion with extrapolated low frequency data, 2016.

\bibitem{vision}
Alex Krizhevsky, Ilya Sutskever, and Geoffrey~E Hinton.
\newblock Imagenet classification with deep convolutional neural networks.
\newblock In F.~Pereira, C.J. Burges, L.~Bottou, and K.Q. Weinberger, editors,
  {\em Advances in Neural Information Processing Systems}, volume~25. Curran
  Associates, Inc., 2012.

\bibitem{language}
Alex Graves.
\newblock Generating sequences with recurrent neural networks, 2013.

\bibitem{speech}
Abdel-rahman Mohamed, George~E. Dahl, and Geoffrey Hinton.
\newblock Acoustic modeling using deep belief networks.
\newblock {\em IEEE Transactions on Audio, Speech, and Language Processing},
  20(1):14--22, 2012.

\bibitem{fno}
Zongyi Li, Nikola Kovachki, Kamyar Azizzadenesheli, Burigede Liu, Kaushik
  Bhattacharya, Andrew Stuart, and Anima Anandkumar.
\newblock Fourier neural operator for parametric partial differential
  equations, 2020.

\bibitem{deeponet}
Lu~Lu, Pengzhan Jin, Guofei Pang, Zhongqiang Zhang, and George~Em Karniadakis.
\newblock Learning nonlinear operators via deeponet based on the universal
  approximation theorem of operators.
\newblock In {\em Advances in Neural Information Processing Systems}, volume~3.
  Nature Machine Intelligence, 2021.

\bibitem{MLZ}
Matthew Li, Laurent Demanet, and Leonardo Zepeda-N\'{u}\~{n}ez.
\newblock Wide-band butterfly network: Stable and efficient inversion via
  multi-frequency neural networks.
\newblock {\em Multiscale Modeling \& Simulation}, 20(4):1191--1227, 2022.

\bibitem{Khoo_YingSwitchNet:2019}
Y.~Khoo and L.~Ying.
\newblock Switch{N}et: A neural network model for forward and inverse
  scattering problems.
\newblock {\em SIAM J. Sci. Comput.}, 41(5):A3182--A3201, 2019.

\bibitem{FanYing:scattering}
Y~Fan and L.~Ying.
\newblock Solving inverse wave scattering with deep learning.
\newblock {\em arXiv:1911.13202}, 2019.

\bibitem{PINN_Inverse_Problems}
Yuyao Chen, Lu~Lu, George~Em Karniadakis, and Luca~Dal Negro.
\newblock Physics-informed neural networks for inverse problems in nano-optics
  and metamaterials.
\newblock {\em Opt. Express}, 28(8):11618--11633, Apr 2020.

\bibitem{MoDL}
H.~K. {Aggarwal}, M.~P. {Mani}, and M.~{Jacob}.
\newblock {MoDL}: Model-based deep learning architecture for inverse problems.
\newblock {\em IEEE Transactions on Medical Imaging}, 38(2):394--405, 2019.

\bibitem{Neumann_Networks}
Davis Gilton, Greg Ongie, and Rebecca Willett.
\newblock Neumann networks for inverse problems in imaging.
\newblock {\em arXiv preprint arXiv:1901.03707}, 2019.

\bibitem{Mao:2016}
Xiaojiao Mao, Chunhua Shen, and Yu-Bin Yang.
\newblock Image restoration using very deep convolutional encoder-decoder
  networks with symmetric skip connections.
\newblock In D.~D. Lee, M.~Sugiyama, U.~V. Luxburg, I.~Guyon, and R.~Garnett,
  editors, {\em Advances in Neural Information Processing Systems 29}, pages
  2802--2810. Curran Associates, Inc., 2016.

\bibitem{U-Net}
Olaf Ronneberger, Philipp Fischer, and Thomas Brox.
\newblock {\em {U-Net}: Convolutional Networks for Biomedical Image
  Segmentation}, pages 234--241.
\newblock Springer International Publishing, Cham, 2015.

\bibitem{Invariant_Scattering_CNN}
J.~{Bruna} and S~{Mallat}.
\newblock Invariant scattering convolution networks.
\newblock {\em IEEE Transactions on Pattern Analysis and Machine Intelligence},
  35(8):1872--1886, 2013.

\bibitem{Framelets}
Jong~Chul Ye, Yoseob Han, and Eunju Cha.
\newblock Deep convolutional framelets: A general deep learning framework for
  inverse problems.
\newblock {\em SIAM Journal on Imaging Sciences}, 11(2):991--1048, 2018.

\bibitem{Framelets_denosing}
E.~{Kang}, W.~{Chang}, J.~{Yoo}, and J.~C. {Ye}.
\newblock Deep convolutional framelet denosing for low-dose {CT} via wavelet
  residual network.
\newblock {\em IEEE Transactions on Medical Imaging}, 37(6):1358--1369, 2018.

\bibitem{MNNH2}
Y.~Fan, J.~Feliu-Fab{\`a}, L.~Lin, L.~Ying, and L.~Zepeda-N{\'u}{\~{n}}ez.
\newblock A multiscale neural network based on hierarchical nested bases.
\newblock {\em Research in the Mathematical Sciences}, 6(2):21, Mar. 2019.

\bibitem{FourierNeuralOp}
Zongyi Li, Nikola~Borislavov Kovachki, Kamyar Azizzadenesheli, Burigede liu,
  Kaushik Bhattacharya, Andrew Stuart, and Anima Anandkumar.
\newblock {Fourier Neural Operator for Parametric Partial Differential
  Equations}.
\newblock In {\em International Conference on Learning Representations}, 2021.

\bibitem{LRC2020}
Yifan Peng, Lin Lin, Lexing Ying, and Leonardo Zepeda-N\'{u}\~{n}ez.
\newblock Efficient long-range convolutions for point clouds, 2020.

\bibitem{linear_butterfly}
Tri Dao, Albert Gu, Matthew Eichhorn, Atri Rudra, and Christopher R{\'e}.
\newblock Learning fast algorithms for linear transforms using butterfly
  factorizations.
\newblock {\em Proceedings of Machine Learning Research}, 97:1517--1527, 06
  2019.

\bibitem{Butterfly-Net2}
Zhongshu Xu, Yingzhou Li, and Xiuyuan Cheng.
\newblock {Butterfly-Net2: Simplified Butterfly-Net and F}ourier transform
  initialization.
\newblock In Jianfeng Lu and Rachel Ward, editors, {\em Proceedings of The
  First Mathematical and Scientific Machine Learning Conference}, volume 107 of
  {\em Proceedings of Machine Learning Research}, pages 431--450, Princeton
  University, Princeton, NJ, USA, 20--24 Jul. 2020. PMLR.

\bibitem{Yingzhou2018}
Y.~Li, X.~Cheng, and J.~Lu.
\newblock Butterfly-{Net}: Optimal function representation based on
  convolutional neural networks.
\newblock {\em arXiv preprint arXiv:1805.07451}, 2018.

\bibitem{Butterfly-factorization:Liu_Xing2020}
Y.~Liu, X.~Xing, H.~Guo, E.~Michielssen, and X.~S. Ghysels, P.~Li.
\newblock Butterfly factorization via randomized matrix-vector multiplications.
\newblock {\em arXiv:2002.03400}, 2020.

\bibitem{optimization}
Carlos Borges, Adrianna Gillman, and Leslie Greengard.
\newblock High resolution inverse scattering in two dimensions using recursive
  linearization, 2016.

\bibitem{Plessix_2006:ajoint_state}
R.-E. Plessix.
\newblock {A review of the adjoint-state method for computing the gradient of a
  functional with geophysical applications}.
\newblock {\em Geophysical Journal International}, 167(2):495--503, 11 2006.

\bibitem{fbp}
D.~Colton and R.~Kress.
\newblock {\em Integral Equation Methods in Scattering Theory}.
\newblock Society for Industrial and Applied Mathematics, Philadelphia, PA,
  2013.

\bibitem{BF}
Yingzhou Li, Haizhao Yang, Eileen~R. Martin, Kenneth~L. Ho, and Lexing Ying.
\newblock Butterfly factorization.
\newblock 13(2):714--732, jan 2015.

\bibitem{cohenc16:EquivariantCNN}
Taco Cohen and Max Welling.
\newblock Group equivariant convolutional networks.
\newblock In Maria~Florina Balcan and Kilian~Q. Weinberger, editors, {\em
  Proceedings of The 33rd International Conference on Machine Learning},
  volume~48 of {\em Proceedings of Machine Learning Research}, pages
  2990--2999, New York, New York, USA, 20--22 Jun 2016. PMLR.

\bibitem{bronstein2021geometric}
Michael~M. Bronstein, Joan Bruna, Taco Cohen, and Petar Veličković.
\newblock Geometric deep learning: Grids, groups, graphs, geodesics, and
  gauges, 2021.

\bibitem{tensorflow2015}
Mart\'{\i}n Abadi, Ashish Agarwal, Paul Barham, Eugene Brevdo, Zhifeng Chen,
  Craig Citro, Greg~S. Corrado, Andy Davis, Jeffrey Dean, Matthieu Devin,
  Sanjay Ghemawat, Ian Goodfellow, Andrew Harp, Geoffrey Irving, Michael Isard,
  Yangqing Jia, Rafal Jozefowicz, Lukasz Kaiser, Manjunath Kudlur, Josh
  Levenberg, Dandelion Man\'{e}, Rajat Monga, Sherry Moore, Derek Murray, Chris
  Olah, Mike Schuster, Jonathon Shlens, Benoit Steiner, Ilya Sutskever, Kunal
  Talwar, Paul Tucker, Vincent Vanhoucke, Vijay Vasudevan, Fernanda Vi\'{e}gas,
  Oriol Vinyals, Pete Warden, Martin Wattenberg, Martin Wicke, Yuan Yu, and
  Xiaoqiang Zheng.
\newblock {TensorFlow}: Large-scale machine learning on heterogeneous systems,
  2015.
\newblock Software available from tensorflow.org.

\bibitem{kingma2015adam}
Diederik Kingma and Jimmy Ba.
\newblock Adam: a method for stochastic optimization.
\newblock In {\em Proceedings of the International Conference on Learning
  Representations (ICLR)}, May 2015.

\bibitem{glorot}
Xavier Glorot and Yoshua Bengio.
\newblock Understanding the difficulty of training deep feedforward neural
  networks.
\newblock In Yee~Whye Teh and Mike Titterington, editors, {\em Proceedings of
  the Thirteenth International Conference on Artificial Intelligence and
  Statistics}, volume~9 of {\em Proceedings of Machine Learning Research},
  pages 249--256, Chia Laguna Resort, Sardinia, Italy, 13--15 May 2010. PMLR.

\bibitem{Berenger:PML}
J.-P. B\'erenger.
\newblock A perfectly matched layer for the absorption of electromagnetic
  waves.
\newblock {\em J. Comput. Phys.}, 114(2):185--200, 1994.

\bibitem{SL}
L.~A. Shepp and B.~F. Logan.
\newblock The fourier reconstruction of a head section.
\newblock {\em IEEE Transactions on Nuclear Science}, 21(3):21--43, 1974.

\bibitem{pytorch}
Adam Paszke, Sam Gross, Francisco Massa, Adam Lerer, James Bradbury, Gregory
  Chanan, Trevor Killeen, Zeming Lin, Natalia Gimelshein, Luca Antiga, Alban
  Desmaison, Andreas Kopf, Edward Yang, Zachary DeVito, Martin Raison, Alykhan
  Tejani, Sasank Chilamkurthy, Benoit Steiner, Lu~Fang, Junjie Bai, and Soumith
  Chintala.
\newblock Pytorch: An imperative style, high-performance deep learning library.
\newblock In {\em Advances in Neural Information Processing Systems 32}, pages
  8024--8035. Curran Associates, Inc., 2019.

\end{thebibliography}

\appendix
\section{Wide-band Butterfly Network}
Wide-band Butterfly Network utilizes the butterfly factorization and Cooley-Tukey FFT algorithm to reduce the number of trainable weights so that it matches the inherent complexity of the problem. The model does not exploit the rotational equivariance property of the problem, and it directly approximates the back-scattering operator by a matrix in the discrete setting.

As we discussed in Section~\ref{sec:propertybutterfly}, the structure of the butterfly factorization is: 
\begin{equation}\label{eqn:bfappendix}
\sfK\approx \sfU^L\sfG^{L-1}\cdots \sfG^{L/2}\sfM^{L/2}(\sfH^{L/2})^*\cdots(\sfH^{L-1})^*(\sfV^L)^*\, ,
\end{equation}
{where $L$ is the level of the butterfly factorization.} When the matrix is applied on the left to the data, the process can be intuitively interpreted as: $\sfV^L$ extracts a local representation of the vector, and then each $\sfH^l$ compresses two adjacent local representations. Upon the application of the switch matrix $\sfM^{L/2}$ that redistributes the representations from the previous step by permuting the vector, each $\sfG^l$ decompress the representation by splitting it into two, which increases the resolution of the representation. Finally, $\sfU^L$ converts the local representations to sampling points. 
In the process, we notice that the resolution of the representation decreases with each $\sfH^l$ being applied. That is where the idea of the Cooley-Tukey FFT algorithm come into effect.  After we compress the data of higher frequency by applying $\sfH^l$ to it, it is natural to merge the obtained representation of lower resolution with the data of lower frequency. The assimilation of the multi-frequency data is done in the process of applying $\sfH^l$. As in our model, the filtering operator is approximated by a 2D convolutional neural network.

{In our experiments, we use the data of the same dimension as was in the original paper~\cite{MLZ}. Hence, according to the paper, we choose the number of levels as 4, the leaf size as 5 and the rank as 3.} {In particular, the level number $L$ determines the number of factors in the butterfly factorization, each of which is approximated by a LocallyConnected2D layer in Tensorflow~\cite{tensorflow2015}.  The leaf size is chosen so that in any $s\times s$ pixels, the data is not oscillatory. Hence, the data can be further compressed and admits a low-rank representation with the desired rank.}

\section{Fourier Neural Operator}
Fourier Neural Operator~\cite{fno} is a data-driven method that learns operators mapping between infinite dimensional spaces. In the model, the {far-field} data is first lifted into a higher dimensional channel space by a neural network, which is usually done by using a shallow fully-connected network. Then, a series of Fourier layers are applied. The Fourier layer is composed of a Fourier transformation, a linear transformation, i.e. multiplication by trainable parameters, to filter out the higher modes, and a inverse Fourier transformation.  Finally, the result is projected back to the target dimension.

{In our experiments, we follow the training approach of the original paper~\cite{fno}. The model lifts the data to 32-dimensional channel space by applying a linear transformation. Then, four Fourier layers are applied. In particular, the number of higher modes chosen in each layer is $12$. Finally it projects the results from the channel space to the target dimension by applying two linear transformations with an activation function placed between them.}


\end{document}